\documentclass[journal]{IEEEtran}
\usepackage{tikz}
\usepackage{graphicx}
\usepackage{dcolumn}
\usepackage{bm}
\usepackage{caption}
\usepackage{comment}
\usepackage{subcaption}
\usepackage{xcolor}
\usepackage{amsfonts}
\usepackage{ulem}
\usepackage{placeins}
\usepackage{amsmath}
\usepackage{cite}
\usepackage{mathabx}
\usepackage{verbatim}
\usepackage{amssymb}

\begin{document}

\title{A Charge Conserving Exponential Predictor Corrector
FEMPIC Formulation for Relativistic Particle Simulations}
\author{Omkar H. Ramachandran,~\IEEEmembership{Graduate Student Member,~IEEE,}
Leo C. Kempel, ~\IEEEmembership{Fellow, ~IEEE,}
John Luginsland, ~\IEEEmembership{Fellow, ~IEEE,}
B. Shanker,~\IEEEmembership{Fellow,~IEEE}
\thanks{ Omkar H. Ramachandran and Leo C. Kempel are with the Department
of Electrical and Computer Engineering, Michigan State University, East Lansing,
MI, 48824.\protect\\
John Luginsland was with MSU's Electrical and Computer Engineering Department.  He is currently with AFRL/Air Force Office of Scientific Research, Arlington, VA 22201.\protect\\
B. Shanker is with the Department of Electrical and Computer Engineering at The Ohio State University, OH, 43210.
E-mail: ramach21@msu.edu}
}
\maketitle

\begin{abstract}
The state of art of charge-conserving electromagnetic finite element particle-in-cell has grown 
by leaps and bounds in the past few years. 
These advances have primarily been achieved for leap-frog time stepping schemes for Maxwell 
solvers, in large part, due to the method strictly following the proper space for representing 
fields, charges, and measuring currents. 
Unfortunately, leap-frog based solvers (and their other incarnations) are only conditionally 
stable. 
Recent advances have made Electromagnetic Finite Element Particle-in-Cell (EM-FEMPIC) methods built around \emph{unconditionally} stable time 
stepping schemes were shown to conserve charge.
Together with the use of a quasi-Helmholtz decomposition, these methods were both unconditionally 
stable and satisfied Gauss' Laws to machine precision. 
However, this architecture was developed for systems with explicit particle integrators where 
fields and velocities were off by a time step. 
While completely self-consistent methods exist in the literature, they follow the classic rubric: 
collect a system of first order differential equations (Maxwell and Newton equations) and use an 
integrator to solve the combined system. 
These methods suffer from the same side-effect as earlier--they are conditionally stable. 
Here we propose a different approach; we pair an unconditionally stable Maxwell solver to an 
exponential predictor-corrector method for Newton's equations. 
As we will show via numerical experiments, the proposed method conserves energy within a PIC scheme, has 
an unconditionally stable EM solve, solves Newton's equations to much higher accuracy than a traditional Boris solver
and conserves charge to machine precision.
We further demonstrate benefits compared to other polynomial methods to solve Newton's equations, 
like the well known Boris push. 
\end{abstract}

\section{Introduction \label{sec:introduction}}

The co-simulation of moving charged particles within an electromagnetic field is of great 
importance in a number of applications, ranging from accelerator design, pulsed power devices,
device sterilization, high precision etching, among many others \cite{marchand2011ptetra,lemke1999three,fourkal2002particle}. This demand has led to the development of several computational methods, foremost among which is the Particle-in-Cell (PIC) technique.
While PIC methods have traditionally based themselves around finite difference time domain (FD-TD)
electromagnetic solvers, recent advances have made it possible to use the finite 
element method (FEM) \cite{pinto2014charge,moon2014exact,squire2012geometric,moon2015exact,o2021set} 
-- and along with it, the ability to better resolve complex geometries and better understanding of the need to strictly follow \emph{de-Rham} \cite{tonti,deschamps} relation for fields, fluxes, charges and currents. These relations rigorously define the correct function spaces that the basis function used to represent various field and current quantities should live in, as well as provide appropriate transformations that allow one to relate different quantities to each other. In the literature, this nuance is also referred to as structure preserving methods \cite{kormann2021energy,xiao_2018,CamposPinto2022}. The advances in Electromagnetic Finite Element Particle-in-Cell methods (EM-FEMPIC) have primarily utilized leapfrog based Maxwell solvers, which can
be shown to be devoid of null-spaces and satisfy Gauss' Laws, thereby 
conserving charge without having to resort to costly post-processing, like divergence cleaning \cite{pinto2014charge}. These properties were thoroughly investigated in \cite{crawford2021rubrics} and a number of 
general rules were provided that, when followed, would result in a method that naturally conserved charge.

As one can imagine, the principal drawback is that methods based on a leap-frog solution of Maxwell's equations are only conditionally stable and the time step size is related to the finest feature in the mesh. It is possible to use \emph{unconditionally} stable implicit Maxwell solvers as well as wave equation solvers for PIC and still conserve charge \cite{Oconnor_time_integration}; however, one needs to handle a DC null space in case of Maxwell solvers or time growing null space for the wave equation. 
However, both can be completely overcome by using a quasi-Helmholtz decomposition for representation of fields and fluxes \cite{oconnor2021quasihelmholtz}. 
Indeed, this approach explicitly enforces \emph{both} of Gauss' laws in terms of the involution equations. 
More importantly,  the solution to the two curl equations have zero divergence \cite{oconnor2021quasihelmholtz} and, therefore, do \emph{not} corrupt the satisfaction of Gauss' law.
This has further been advanced to utilize techniques from envelope tracking to better analyse narrow band systems \cite{ramachandran2022envelopetracking}.

Despite these advances, there remain some unsolved problems. 
One of the principal challenges among these is to ensure self-consistency at every time step; i.e the fields, positions and velocities need to be self consistent. 
What is typically done is to reduce both Maxwell and Newton equations to a set of coupled first order partial differential equations. 
These are then collectively evolved using an appropriate integrator. Alternately, one offsets the evolution of Maxwell's equations from that of Newton's by one time step, i.e., the evolution is explicit. 
While the former is self-consistent, the latter is not. Delving into the latter for the sake of exposition, we note the following: (a) evolution of Newton's laws using Boris push (a very robust and popular technique) and leap-frog for Maxwell systems is one of the most popular approaches in setting up a PIC scheme \cite{pinto2014charge}; (b) Boris time stepping is either momentum or energy conserving, preserves the volume of phase space in either case,  and has low computational overhead \cite{boris,why_boris_is_good}; 
(c) phase errors are known to exist as one increases the timestep size \cite{2018_boris_alt_method};  (d) leap-frog evolution for Maxwell systems is \emph{conditionally stable} with time step sizes depending on the smallest feature. 
There is significant interest in constructing methods that are self-consistent \cite{Bacchini_2019}.
Existing implementations of such methods typically follow a so-called `kinetic enslavement' 
approach, with Maxwell and particle systems combined together into a large coupled solver,
which can then be evaluated through an appropriate nonlinear solver \cite{markidis_2011}.
Such an approach comes with the obvious downside of being very computationally expensive, since
the size of the system is in the order of the number of particles.
As a result, many implementations of kinetic enslavement often rely on matrix-free Jacobian
free nonlinear solvers to avoid the enormous storage footprint of the combined system matrix 
\cite{implicit_thesis}. 
As is evident, relying on a specific nonlinear solver comes with downsides. Further, time step sizes used for kinetic enslavement is also subject to stability constraints, and much like the Boris push, position and velocity updates are polynomial based with convergence rates following the order of approximation.


The primary goal of this work is to investigate another approach for self-consistent evolution of particle and Maxwell equations with the following requirements; (a) the evolution of Maxwell systems should be \emph{unconditionally} stable, i.e., the time step sizes should really be only governed by physics and not feature size used to describe the geometry; (b) it should naturally satisfy Gauss' Laws, thereby conserving charge and (c) it should show superior error convergence.

The objectives of this work, therefore, are twofold: (a) First, we analyse the accuracy and error convergence for a novel predictor-corrector integrator based on exponential basis sets \cite{Glaser2009} over a number of particle-only numerical experiments.
(b) Using this integrator, we construct a charge conserving EM-FEMPIC scheme that functions under
a predictor corrector update.
We validate both the proposed integrator and the overall PIC scheme through a number of analytical
and numerical results.


The outline of the paper will be as follows: Sec. \ref{sec:discrete} and Sec. \ref{sec:quasihelmholtz} 
will outline the construction of the quasi-Helmholtz based EM solver. 
Sec. \ref{sec:particle_evol} and following will detail and analyze energy conservation for a 
number of different particle evolution schemes. Sec. \ref{sec:results} will then validate both
the proposed particle evolution scheme and the overall PIC method through a number of analytical
and numerical tests, before concluding remarks are presented in Sec. \ref{sec:summary}.

\section{Formulation \label{sec:formulation}}

\subsection{Problem Statement}

Consider a region $\Omega \in \mathbb{R}^{3}$ bounded by a surface $\partial \Omega$ containing a single charged species. This region is subjected to an external field due to which the charged species accelerate, and in turn produce spatially and temporally varying electric and magnetic fields denoted by $\mathbf{E}(\mathbf{r},t)$ and $\mathbf{B}(\mathbf{r},t)$, respectively, with $\mathbf{r}\in\Omega$ and $t\in[0,\infty)$. The dynamics of the particles in phase space can be represented by a distribution function (PSDF) $f(t,\mathbf{r},\mathbf{v})$ that follows the Vlasov equation:
\begin{equation}\label{eq:distFn}
\begin{split}
    \partial_{t} f(t,\mathbf{r},\mathbf{v}) +\mathbf{v} \cdot \nabla f(t,\mathbf{r},\mathbf{v}) +& \\
    \frac{q}{m}\left[\mathbf{E}(\mathbf{r},t)+\mathbf{v}\times\mathbf{B}(\mathbf{r},t)\right]\cdot \nabla_{v}f(t,\mathbf{r},\mathbf{v}) =& 0.
\end{split}
\end{equation}
$q$ and $m$ here refer to the particle charge and mass respectively. In what follows, we assume that the background media in $\Omega$ is free space. As a result, we denote the permittivity and permeability of free space by $\epsilon_{0}$ and $\mu_{0}$, respectively, and the speed of light by $c$.
\subsection{Discrete Solution \label{sec:discrete}}
Solutions to \eqref{eq:distFn} is found through a Particle-in-Cell formulation, wherein the charge and current distributions within the domain are represented in terms of a finite number of moments of the PSDF $f(t,\mathbf{v},\mathbf{r})$
\begin{subequations}
\begin{equation}
    \rho(\mathbf{r},t) = \sum_{p=1}^{N_{p}} q_{p} S(\mathbf{r}-\mathbf{r}_{p}(t)),
\end{equation}
\begin{equation}
    \mathbf{J}(\mathbf{r},t) = \sum_{p=1}^{N_{p}} q_{p} \mathbf{v}_{p}(t) S(\mathbf{r}-\mathbf{r}_{p}(t)).
\end{equation}
\end{subequations}
Here, $N_{p}$ refers to the total number of such samples or macroparticles in the simulation. And likewise $q_{p}$ and $m_{p}$ refer to the charge and mass of the $p$th particle respectively.
The dynamics of the fields within the simulation domain $\Omega$ obey Maxwell's Equations which can be written as
\begin{subequations}
    \label{eq:maxwell_continuous}
    \begin{equation}
        \nabla \times \mathbf{E}(\mathbf{r},t) = -\partial_{t} \mathbf{B}(\mathbf{r},t),
    \end{equation}
    \begin{equation}
        \nabla \times \mu_{0}^{-1}\mathbf{B}(\mathbf{r},t) = \partial_{t} \mathbf{G}(\mathbf{r},t) + \partial_{t}\epsilon_{0} \mathbf{E}(\mathbf{r},t),
    \end{equation}
\end{subequations}
where $\mathbf{G}(\mathbf{r},t)=\intop_{0}^{t}\mathbf{J}(\mathbf{r},\tau) d\tau$.
Further, the fields are constrained by Gauss' electric and magnetic laws
\begin{subequations}
    \label{eq:gauss_continuous}
    \begin{equation}
        \nabla \cdot \epsilon_{0} \mathbf{E}(\mathbf{r},t) = \rho(\mathbf{r},t),
    \end{equation}
    \begin{equation}
        \nabla \cdot \mathbf{B}(\mathbf{r},t) = 0.
    \end{equation}
\end{subequations}
The curl equations in \eqref{eq:maxwell_continuous} are discretized in space using a finite element formulation over a tetrahedral mesh with $N_{e}$ edges and $N_{f}$ faces, 
wherein the fields are represented in space using Whitney basis forms, $\mathbf{E}(\mathbf{r},t)=\sum_{i=1}^{N_{e}} e_{i}(t) \mathbf{W}^{(1)}_{i}(\mathbf{r})$, $\mathbf{B}(\mathbf{r},t) = \sum_{i=1}^{N_{f}} b_{i}(t) \mathbf{W}^{(2)}_{i}(\mathbf{r})$, for some time-varying coefficient $e_{i}(t)$ and $b_{i}(t)$ respectively and subsequently Galerkin tested. 
Temporally, this system is discretized using a Newmark-$\beta$ \cite{zienkiewicz1977new} formulation with $\gamma=0.5$ and $\beta=0.25$. In what follows, $\bar{E} = \left [e_1 (t), e_2(t), \cdots,  e_{N_e}(t)\right ] $, $\bar{B} = \left [b_1 (t), b_2(t), \cdots,  b_{N_f}(t)\right ] $. We denote evaluations of quantities at $t^n = n\Delta_t$ where $\Delta_t$ is the time step size, via a superscript $n$ such as $\bar{E}^{n}=\bar{E}(t^{n})$.
Complete details on the mixed FEM discretization can be found in \cite{Boss88,crawford2020unconditionally,Wong95,He06,Kett99} and the references therein.

\subsection{Quasi Helmholtz \label{sec:quasihelmholtz}}

Using the formulation described in Sec. \ref{sec:discrete} with an implicit time-marching scheme is well known to not satisfy Gauss' Laws due to null spaces for both Maxwell solver or the wave-equation formulation \cite{wang2010application,venkatarayalu2006removal}, both with implicit unconditionally stable time stepping. 
The most comprehensive method to mitigate this issue is to use a so-called Quasi-Helmholtz decomposition to partition the field components into solenoidal and non-solenoidal components. The non-solenoidal components can then be explicitly solved for, thereby ensuring exact satisfaction of the Coulomb gauge.  What follows is a brief overview of the quasi-Helmholtz setup described in great detail in \cite{oconnor2021quasihelmholtz}.

In what follows, variables with the subscript `$ns$' will refer to non-solenoidal quantities and `$s$' will likewise denote solenoidal quantities.
\subsubsection{Relevant Matrix Definitions\label{sec:qhmats}}
To begin, the sets $\mathcal{N}$, $\mathcal{E}$, $\mathcal{F}$ and $\mathcal{T}$ are defined as the set of nodes, edges, faces and tets within the mesh of the simulation domain respectively having $N_{n}$, $N_{f}$, $N_{e}$ and $N_{t}$ elements.
The various submatrices that will be used in describing the Quasi-Helmholtz framework are as follows:
\begin{subequations}
\begin{equation}
    [\star_\epsilon]_{i,j} = \langle \mathbf{W}^{(1)}_i(\mathbf{r}),\varepsilon\cdot\mathbf{W}^{(1)}_j(\mathbf{r}) \rangle;  i,j \in \mathcal{E}
\end{equation}
\begin{equation}
    [\star_{\mu^{-1}}]_{i,j} = \langle \mathbf{W}^{(2)}_i(\mathbf{r}),\mu^{-1}\cdot\mathbf{W}^{(2)}_j(\mathbf{r})\rangle; i,j \in \mathcal{F}
\end{equation}
\begin{equation}
    [\star_{\rho}]_{i,j} = \langle W^{(3)}_i(\mathbf{r}),W^{(3)}_j(\mathbf{r})\rangle; i,j \in \mathcal{T}
\end{equation}
\end{subequations}
where $\mathbf{W}^{(1)}_{i}$, $\mathbf{W}^{(2)}_{i}$ and $W^{(3)}_{i}$ are the Whitney edge, face and volume basis functions respectively. 
Further, we define the following matrices:
\begin{subequations}
\begin{equation}
    [\mathbf{M}_g]_{i,j} = \langle \mathbf{W}^{(1)}_i (\mathbf{r}) , \nabla W^{(0)}_j (\mathbf{r})\rangle; i \in \mathcal{E}, j \in \mathcal{N} \label{eq:mg}
    \end{equation}
    \begin{equation}
        [\mathbf{M}_c]_{i,j} = \langle  \mathbf{W}^{(2)}_i (\mathbf{r}), [\mathbf{\nabla\times}] \mathbf{W}^{(1)}_j (\mathbf{r})\rangle;  i \in \mathcal{F}, j \in \mathcal{E}
        \label{eq:mc}
    \end{equation}
    \begin{equation}
        [\mathbf{M}_d]_{i,j} = \langle W^{(3)}_i(\mathbf{r}) , [\mathbf{\nabla} \cdot] \mathbf{W}^{(2)}_j (\mathbf{r}) \rangle;  i \in \mathcal{T}, j \in \mathcal{F} \label{eq:md}
    \end{equation}
    \begin{equation}
        [\mathbf{\nabla} ] = \varepsilon[\star_{\epsilon}]^{-1}[\mathbf{M}_g] \label{eq:grad_mat}
    \end{equation}
    \begin{equation}
        [\mathbf{\nabla\times} ] = \mu^{-1}[\star_{\mu^{-1}}]^{-1}[\mathbf{M}_c]\label{eq:curl_mat}
    \end{equation}
    \begin{equation}
        [\mathbf{\nabla}\cdot ]  = [\star_{\rho}]^{-1}[\mathbf{M}_d] \label{eq:div_mat}
    \end{equation}
\end{subequations}

\subsubsection{Projector Definitions}
To separate the non-solenoidal components from the basis forms used for representing the electric field, we define projectors $[\bar{P}]^{\Sigma}_{e}$ and $[\bar{P}]^{\Lambda}_{e}$ to break up the electric field as follows (where $\dagger$ represents a Moore-Penrose pseudoinverse):
\begin{subequations}
\begin{align}
    [\bar{P}]_e^{\Sigma} & = \Sigma(\Sigma^T \Sigma)^{\dagger} \Sigma^{T} \label{eq:proj_es}\\
    [\bar{P}]^{\Lambda}_e & = \mathcal{I} - [\bar{P}]_e^{\Sigma} \label{eq:proj_el}
\end{align}
\begin{equation}
\label{eq:bfld_projector}
\left [ \bar{P}\right ]^\Lambda_b = \mathcal{I} - \Sigma_m\left ( \Sigma_m^T \Sigma_m \right )^{\dagger} \Sigma_m^T 
\end{equation}
\end{subequations}
where $[\Sigma]=\epsilon_{0}[\mathbf{M}_{g}]$ and $[\Sigma]_m=[\nabla \cdot]^T$. 
Using these projectors, we can now define a decomposition for the electric flux density as
\begin{equation}
\label{eq:efld_projected}
\bar{D}(t) = \Sigma \bar{E}_{ns}(t) + \left [ \bar{P}\right ] ^\Lambda_e \bar{D}(t)
\end{equation}
and the magnetic flux density as
\begin{equation}
\label{eq:bfld_projected}
\bar{B}_{s}(t) = \left [ \bar{P}\right ] ^\Lambda_b \bar{B}(t)
\end{equation}
Put simply, these projectors allow us to separate out $\bar{B}(t)$ and $\bar{D}(t)$ into terms that have exactly zero divergence (i.e, the 'solenoidal components' $\bar{B}_{s}(t)$ and $\left [ \bar{P}\right ] ^\Lambda_e \bar{D}(t)$) and terms with non-zero divergence ($\Sigma \bar{E}_{ns}(t)$).
Since the application of a discrete divergence operator on either projector is 0 (as demonstrated in \cite{oconnor2021quasihelmholtz}), one can see that the formulation forces the divergence of $\bar{B}_s(t)$ to zero.

\subsubsection{Discrete System}

As stated earlier, our ultimate goal is to solve \emph{all} of Maxwell's equations. First given the decomposition, the application of the discrete divergence on \eqref{eq:bfld_projected} and \eqref{eq:efld_projected} results in an identically zero matrix and
\begin{equation}
\label{eq:ns_equation}
\left[C_{z}^{e}\right]^{T}\left[\nabla\right]^{T}\left[\star_{e}\right]\left[\nabla\right]\left[C_{z}^{e}\right]\bar{E}_{ns}(t)
=-\left[C_{z}^{e}\right]^{T}\left[\nabla\right]^{T}\bar{G}(t).
\end{equation}
Upon rewriting Maxwell's equations to make use of the decomposition in \eqref{eq:efld_projected} and \eqref{eq:bfld_projected}; and with $\bar{E}_{ns}(t)$ in \eqref{eq:ns_equation} explicitly solved for, we obtain 
\begin{equation}
    \label{eq:sol_disc}
    \begin{split}
        \left[\bar{Z}\right]_{11}\partial_{t}\bar{B}_{s}(t)
        +&\left[\bar{Z}\right]_{12}\bar{E}_{s}(t)=
        -\left[\bar{Z}\right]_{13}\bar{E}_{ns}(t)\\
        \left[\bar{Z}\right]_{21}\partial_{t}\bar{E}_{s}(t)
        - &
        \left[\bar{Z}\right]_{22}\bar{B}_{s}(t)
        =-\partial_{t}\bar{G}(t)
        -\left[\bar{Z}\right]_{23}\partial_{t}\bar{E}_{ns}(t)
    \end{split}
\end{equation}
where $\bar{E}_{s}$ and $\bar{B}_{s}$ are extracted from the simply connected mesh using tree-cotree maps $\left[C^{e}_{c}\right]$ and $\left[C^{b}_{c}\right]$. Constructing this mapping is trivial for simply connected structures, but is trickier for multiply connected geometries. The various submatrices involved in \eqref{eq:sol_disc} are defined as:
\begin{subequations}
\begin{equation}
    [\mathbf{Z}]_{11} =   [\mathbf{C}_c^b]^T[\mathbf{P}]_b^\Lambda [\mathbf{C}_c^b]
    \end{equation}
    \begin{equation}
    [\mathbf{Z}]_{12} =  [\mathbf{C}_c^b]^T [\mathbf{\nabla\times}] [\star_\varepsilon]^{-1}[\mathbf{P}]_e^\Lambda [\star_\varepsilon] [\mathbf{C}_c^e]
    \end{equation}
    \begin{equation}
    [\mathbf{Z}]_{13} =  [\mathbf{C}_c^b]^T [\mathbf{\nabla\times}] [\star_\varepsilon]^{-1}\Sigma [\mathbf{C}_z^e]
    \end{equation}
    \begin{equation}
    [\mathbf{Z}]_{21}  =  [\mathbf{C}_c^e]^{T}[\mathbf{P}]_e^\Lambda[\star_\varepsilon][\mathbf{C}_c^e]
    \end{equation}
    \begin{equation}
    [\mathbf{Z}]_{22} =  [\mathbf{C}_c^e]^{T} [\mathbf{\nabla\times}]^T [\star_{\mu^{-1}}][\mathbf{P}]_b^\Lambda [\mathbf{C}_c^b]
    \end{equation}
    \begin{equation}
    [\mathbf{Z}]_{23} =  [\mathbf{C}_c^e]^{T} \Sigma  [\mathbf{C}_z^e]
\end{equation}
\end{subequations}
As alluded to earlier, the solution to  \eqref{eq:sol_disc} is obtained using a Newmark-$\beta$ time stepping scheme. The solution to $\bar{E}_{s} (t)$ and $\bar{B}_{s}  (t)$ still has a $DC$ null space, but their divergence is exactly zero. As a result, Gauss' laws are still exactly satisfied. 

\subsection{Particle Evolution \label{sec:particle_evol}}

A particle evolution scheme is primarily involved in obtaining solutions in space and time
of the position and velocity of a particle in response to a prescribed set of electromagnetic
fields. Specifically, the particle motion can be obtained by
\begin{subequations}
\begin{equation}
    \frac{\partial \gamma[\mathbf{v}_{p}]\mathbf{v}_{p} (t)}{\partial t} = \mathbf{a}_{p}(t) = \frac{q}{m} 
        \left(\mathbf{E}(\mathbf{r}_{p},t) + \mathbf{v}_{p}\times\mathbf{B}(\mathbf{r}_{p},t)\right)
    \label{eq:lorentz_v}
\end{equation}
\begin{equation}
    \frac{\partial \mathbf{r}_{p}(t)}{\partial t} = \mathbf{v}_{p}(t)
    \label{eq:lorentz_r}
\end{equation}
    \label{eq:lorentz}
\end{subequations}
where $\mathbf{v}_{p}(t)$ and $\mathbf{r}_{p}(t)$ are each vectors in $\mathbb{R}^{3}$ and 
refer to a given particle's velocity and position respectively at time $t$, and  
$\gamma[\mathbf{v}]=(1-\left|\mathbf{v}\right|^{2}/c^{2})$ is the relativistic time dialation
functional. In what follows, we analyze different approaches to solving these equations as well as their stability properties. 

\subsection{Boris Push}
The Boris method is fundamentally a modified leap-frog solver for \eqref{eq:lorentz_v} and \eqref{eq:lorentz_r} that breaks the update for the relativistic velocity $\mathbf{u}_{p}$ into three fundamental steps.
Defining $\mathbf{u}_{p}^{n}$ and $\mathbf{u}_{p}^{n+1}$ as particle 3-velocities at timestep $n$,
\begin{enumerate}
    \item First, half of the electric force is added to a dummy variable $\mathbf{u}^{-}_{p}$.
        \begin{equation}
            \label{eq:first_half}
            \mathbf{u}^{-}_{p} = \mathbf{u}_{p}^{n} + \mathbf{\epsilon}\Delta_{t}
        \end{equation}
        where $\mathbf{\epsilon} = (q_{p}/2m_{p})\mathbf{E}(\mathbf{r}_{p}^{n+1/2},t^{n+1/2})$ and $\Delta_{t}$ is the timestep size.
    \item Next, the magnetic contribution is added,
        \begin{equation}
            \mathbf{u}^{+}_{p} = \mathbf{u}^{-}_{p} + \Delta_{t}\frac{q_{p}}{m_{p}}\left(\mathbf{v}_{p}^{n+1/2}\times \mathbf{B}(\mathbf{r}_{p}^{n+1/2},t^{n+1/2})\right)
            \label{eq:rotation}
        \end{equation}
    \item Finally, the other half of the electric force is added to obtain $\mathbf{u}^{n+1}_{p}$.
        \begin{equation}
            \label{eq:final_half}
            \mathbf{u}_{p}^{n+1} = \mathbf{u}^{+}_{p} + \mathbf{\epsilon}\Delta_{t}
        \end{equation}
\end{enumerate}
While different ways of evaluating \eqref{eq:rotation} have been proposed over the years, with some introducing phase errors \cite{2018_boris_alt_method}, they all share the fundamental property that they execute a
perfect rotation of $\mathbf{u}^{-}_{p}$ and as a result, an applied magnetic force will never
do work on a particle.
Hence, the method as a whole is energy conserving under the action of an arbitrary magnetic field.
While many different methods exist to prove energy conservation, we considered the $z$-transform of the single step update for the implicit leap-frog representation with the forcing terms evaluated at half steps.
Proving that this stencil conserves energy is sufficient to show that \eqref{eq:first_half}-\eqref{eq:final_half} broken up into substeps will conserve energy.
We can write the velocity update stencil for a particle moving due to a general magnetic field as follows:
\begin{equation}
    \left(\mathbf{I}-\mathbf{\Omega}\right)\mathbf{v}^{n+1}_{p} = \left(\mathbf{I}+\mathbf{\Omega}\right)\mathbf{v}^{n}_{p}.
    \label{eq:update_boris_B}
\end{equation}
Here $\mathbf{\Omega}$ can be written as
\begin{equation}
    \label{eq:omega_w}
    \mathbf{\Omega} = \left(\begin{matrix}
            0 & \omega_{z} & -\omega_{y} \\
            -\omega_{z} & 0 & \omega_{x} \\
            \omega_{y} & -\omega_{x} & 0
        \end{matrix}\right),
\end{equation}
where $\omega_{k} = \Delta_{t} (q/m\gamma_{p}^{n+1/2})\mathbf{B}_{k}(\mathbf{r}_{p}^{n+1/2},t^{n+1/2})\;\;\text{for}\; k\in{x,y,z}$ and $\mathbf{I}$ is the regular identity matrix in $\mathbb{R}^{3\times 3}$.
Defining $\mathbf{M}=\left(\mathbf{I}+\mathbf{\Omega}\right)^{-1}\left(\mathbf{I}-\mathbf{\Omega}\right)$, the $z$-transform of \eqref{eq:update_boris_B} becomes
\begin{equation}
    \mathbf{M} z= \lambda z = 1
    \label{eq:Boris_Z_roots}
\end{equation}
where $\lambda$ refers to the eigenvalues of $\mathbf{M}$. 
Under the spectral theorem, for a system to be stable, $|z| <= 1$ for any combination of eigenvalues chosen. However, we can define a stronger condition
$|z| = 1$ wherein the update stencil strictly defines a perfect rotation on the complex plane under the action of a magnetic field. 
This approach has been used in the past to appraise numerical methods for energy conservation; see \cite{zienkiewicz1977new,crawford2020unconditionally} and the references therein.
As a result, $\left|\left(\mathbf{v}_{p}^{n+1}\right)^{2}\right| = \left|\left(\mathbf{v}_{p}^{n}\right)^{2}\right|$ which implies that the particle's kinetic energy is exactly conserved regardless of the form of the magnetic field.
For the Boris update, the eigenvalues $\lambda$ are 
\begin{equation}
    \lambda = \left\{\begin{matrix}
            1 \\
            \frac{1 - \left(\omega_{x}^{2} + \omega_{x}^{2} + \omega_{x}^{2}\right) - 2j\sqrt{\omega_{x}^{2} + \omega_{x}^{2} + \omega_{x}^{2}}}{1+\omega_{x}^{2} + \omega_{y}^{2} + \omega_{z}^{2}} \\
            \frac{1 - \left(\omega_{x}^{2} + \omega_{x}^{2} + \omega_{x}^{2}\right) + 2j\sqrt{\omega_{x}^{2} + \omega_{x}^{2} + \omega_{x}^{2}}}{1+\omega_{x}^{2} + \omega_{y}^{2} + \omega_{z}^{2}}
    \end{matrix}\right\}
    \label{eq:eigenvalue_sols_Boris}
\end{equation}
It can be shown through simple algebraic manipulation that that each of the three eigenvalues in \eqref{eq:eigenvalue_sols_Boris} has exactly unit magnitude for any $\omega_{x},\omega_{y},\omega_{z} \in \mathbb{R}$.
Therefore, when combined with \eqref{eq:Boris_Z_roots}, we see that the roots of the $z$-transform has to lie exactly on the unit circle, thereby conserving energy.
\subsection{Polynomial Predictor Corrector}
Next, we describe the Adams predictor-corrector update. 
These methods are generally derived by representing the forcing terms $\mathbf{a}_{p}(t)=q_{p}/m_{p}\left(\mathbf{E}(\mathbf{r}_{p},t) + \mathbf{v}_{p}(t)\times\mathbf{B}(\mathbf{r}_{p},t)\right)$ and $\mathbf{v}_{p}(t)$
in terms of polynomial functions for a prescribed order and testing them with a pulse function defined over the length of one timestep.
Specifically, the explicit predictor step of the process (typically referred to as an Adams-Bashforth update), involves writing the acceleration $\mathbf{a}_{p}$ and relativistic velocity $\mathbf{u}_{p}$ as
\begin{subequations}
    \label{eq:representation}
    \begin{equation}
        \mathbf{a}_{p}(t) = \sum_{i=1}^{N_{t}}\mathbf{a}^{i}_{p} L_{k}(t-t^{i})
    \end{equation}
    \begin{equation}
        \mathbf{v}_{p}(t) = \sum_{i=1}^{N_{t}}\mathbf{v}_{p}^{i} L_{k}(t-t^{i})
    \end{equation}
\end{subequations}
where $L_{k}(t-t_{i})$ refers to a $k$th order Lagrange polynomial involving the sample at $t_{i}$ and the $k-1$ samples before it in time.
Next, a stencil is obtained from the differential equation by taking an inner product with a  pulse function
\begin{equation}
    P^{n}(t) = \begin{cases}
        1 & \text{if}\; t\in\left[t^{n},t^{n+1}\right] \\
        0 & \text{otherwise}
    \end{cases},
\end{equation}
giving us the following update stencil upon evaluation of the requisite integrals for $k=4$:
\begin{subequations}\label{eq:adam-bashforth}
    \begin{equation}
    \mathbf{u}_{p}^{n+1} = \mathbf{u}_{p}^{n} + \frac{\Delta_t}{24}(55 \mathbf{a}_p^n - 59 \mathbf{a}_p^{n-1} + 37 \mathbf{a}_{p}^{n-2}-9\mathbf{a}_{p}^{n-3}),
    \end{equation}
    \begin{equation}
        \mathbf{r}_{p}^{n+1} = \mathbf{r}_{p}^{n} + \frac{\Delta_t}{24}(55 \mathbf{v}_{p}^n - 59 \mathbf{v}_{p}^{n-1} + 37 \mathbf{v}_{p}^{n-2}-9\mathbf{v}_{p}^{n-3}).
        \label{eq:adam-bashforth-pos}
    \end{equation}
\end{subequations}
One can similarly device an \emph{implicit} scheme by modifying \eqref{eq:representation} such that the Lagrange polynomials center on $n+1$ and go back $k-1$ steps. 
Upon using the same testing function, we obtain a so-called Adams-Moulton update, which for $k=4$ looks as follows:
\begin{subequations}\label{eq:adams-moulton}
    \begin{equation}
        \mathbf{u}_{p}^{n+1} = \mathbf{u}_{p}^{n} + \frac{\Delta_t}{24}(9 \mathbf{a}_p^{n+1} + 19 \mathbf{a}_p^{n} - 5 \mathbf{a}_{p}^{n-1} + \mathbf{a}_{p}^{n-2})
    \end{equation}
    \begin{equation}
        \mathbf{r}_{p}^{n+1} = \mathbf{r}_{p}^{n} + \frac{\Delta_t}{24}(9 \mathbf{v}_{p}^{n+1} + 19 \mathbf{v}_{p}^{n} - 9 \mathbf{v}_{p}^{n-1} + \mathbf{v}_{p}^{n-2}).
        \label{eq:adam-moulton-pos}
    \end{equation}
\end{subequations}
To construct a predictor-corrector scheme using these stencils, at each timestep, a guess velocity $\tilde{\mathbf{v}}_{p}$ and position $\tilde{\mathbf{r}}_{p}$ are computed following \eqref{eq:adam-bashforth}.
These are then inserted into \eqref{eq:adams-moulton} in place of $\mathbf{v}^{n+1}_{p}$ and $\mathbf{r}^{n+1}_{p}$. Going forward, we shall refer to this predictor corrector setup with fourth order polynomials as 'Adams4', and likewise refer to a similar method using third order polynomials as 'Adams3'.

\subsubsection{Analysis of Energy Conservation}

We apply the same means of analysing energy conservation here as we did for the Boris update. 
Consider a simple system wherein the magnetic field is constant and only pointing along the $\hat{z}$ direction.
To give the method the best possible chance of conserving energy, let us assume that $\mathbf{u}_{p}^{n}$, $\mathbf{u}_{p}^{n-1}$ and $\mathbf{u}_{p}^{n-2}$ have identical magnitude.
Therefore, $\gamma_{p}^{n}=\gamma_{p}^{n-1}=\gamma_{p}^{n-2}$.
If under these very favorable conditions, the method fails to conserve energy, then it is sufficient to assume that it will fail to do so under a general field.
Under our prescriptions, the corrector step of the Adams4 velocity update can be written as follows
\begin{equation}
    \left(\mathbf{I}+9\mathbf{\Omega}\right)\mathbf{u}_{p}^{n+1} - \left(\mathbf{I}-19\mathbf{\Omega}\right)\mathbf{u}_{p}^{n} + 5\mathbf{\Omega}\mathbf{u}_{p}^{n-1} - \mathbf{\Omega}\mathbf{u}_{p}^{n-2} = 0 
\end{equation}
where $\mathbf{\Omega}$ has the same form as \eqref{eq:omega_w}, but with $\omega_{x}=\omega_{y}=0$ and $\omega_{z}=(\Delta_{t}q_{p}/24m_{p})\gamma_{p}^{n}\mathbf{B}(\mathbf{r}_{p}^{n},t^{n})$.
Replacing each matrix product with its eigenvalues and taking a $z$-transform gives us
\begin{equation}
    \left(\lambda_{1}z^{3}-\lambda_{2}z^{2}+5\lambda_{3}z+\lambda_{4}\right)\mathbf{u}_{p}(z)=0
\end{equation}
For our chosen system, the eigenvalues can take the following values:
\begin{subequations}
    \begin{equation}
        \lambda_{1} = \left\{1, 1-j9\omega_{z},1+j9\omega_{z}\right\}
    \end{equation}
    \begin{equation}
        \lambda_{2} = \left\{1, 1+j19\omega_{z},1-j19\omega_{z}\right\}
    \end{equation}
    \begin{equation}
        \lambda_{3} = \lambda_{4} = \left\{0, j\omega_{z},-j\omega_{z}\right\}
    \end{equation}
\end{subequations}
By inspection, it is fairly clear that only one combination of eigenvalue choices yields a root for $z$ that is on the unit circle, i.e when $\lambda_{1}=\lambda_{2}=1$ and $\lambda_{3}=\lambda_{4}=0$, giving us,
\begin{equation}
    z^{3}-z^{2} = z^{2}(z-1)=0
\end{equation}
which admits a root of $z=1$.
Every other combination of eigenvalues leads to roots with magnitudes larger or smaller than 1 in the general case.
While one might be able to find specific values of $\omega_{z}$ for which multiple combinations of eigenvalues yield roots on the unit sphere, this would make the method as a whole extremely sensitive to step-size and generally impractical as an energy-conserving method to advance particles.
\subsection{Exponential Predictor-Corrector}
So far, we have considered update schemes that use polynomials to interpolate the particle parameter samples.
However, it is possible to use exponential functions for this purpose. 
One method to construct a $PE(CE)^{m}$ method for ODEs was proposed in \cite{Glaser2009}, wherein the update for the relativistic velocity and position expressions can be written as
\begin{subequations}
    \label{eq:update}
\begin{equation}
    \mathbf{u}^{n+1}_{p} = \sum_{i=1}^{k} \left[p_{i} \mathbf{v}^{n+1-i}_{p} + p_{k+i} \mathbf{a}_{p}^{n+1-i}\right]
\end{equation}
\begin{equation}
    \mathbf{r}_{p}^{n+1} = \sum_{i=1}^{k} \left[p_{i} \mathbf{r}_{p}^{n+1-i} + p_{k+i} \mathbf{v}_{p}^{n+1-i}\right]
\end{equation}
\end{subequations}
for an appropriate set of coefficients $p_{i}$, the specifics of constructing which can be found in \cite{Glaser2009}.
This choice of coefficients encodes an exponential interpolation for the unknowns $\mathbf{v}_{p}(t)$ and $\mathbf{r}_{p}(t)$ in the interval $[t^{n},t^{n+k}]$, which can be written as
\begin{subequations}
    \begin{equation}
        \mathbf{u}_{p}(t) = \sum_{i=n}^{n+k} \alpha^{u}_{i}\mathbf{u}_{p}^{i} e^{\lambda_{i} t}
    \end{equation}
    \begin{equation}
        \mathbf{r}_{p}(t) = \sum_{i=n}^{n+k} \alpha^{r}_{i}\mathbf{r}_{p}^{i} e^{\lambda_{i} t}
    \end{equation}
\end{subequations}
where $\alpha_{k}, \lambda_{k} \in \mathbb{C}$ refer to appropriately chosen coefficients in the complex plane for the parameters in the setup. 
Given the number of prior timesteps involved at any given update, it is not possible to prove stability/energy conservation for an arbitrary applied field. However, unlike the polynomial predictor corrector method, we show through numerical examples that the method does conserve energy for both single particle systems and when incorporated into a full EM-FEMPIC simulation. 
\subsubsection{Computational Complexity}
Since the coefficients for the exponential PC method are not known analytically, there is computational overhead in computing them. However, solving for $p_{i}$ only needs to be done once after which it can be tabulated and used to solve different ODEs discretized using the same number of basis sets.
The costs, therefore, of using this scheme within a PIC method is identical to any $m$-step predictor corrector method. Suppose one uses an iterative solver for the fields, the complexity of evolving Maxwell's equations at any given timestep can be written as $\mathcal{O}(N_{\text{iter}}N_{\text{dof}})$ where $N_{\text{iter}}$ denotes the number of iterations the solver takes to converge and $N_{\text{dof}}$ the total number of unknowns in the system. 
In practice, this cost is significantly higher than the particle evolution. Furthermore, the unconditionally stable nature of the field solve allows for taking much larger timesteps than would otherwise be possible using an explicit scheme.
As a result, even though the Boris and Adams particle update require fewer multiplications than the PC1 scheme used in the paper, the added complexity is negligible in relation to the cost associated with solving for the fields.

\section{Results \label{sec:results}}
The numerical examples in this section are broken up into three parts. (1) First, we consider systems forced by uniform electric and magnetic fields; (2) We then consider systems with nonuniform fields; (3) Finally, the entire EM-FEMPIC scheme is validated through an expanding particle beam simulation.
\subsection{Uniform fields}
In this section, we will consider the motion of particles in response to a combination of electric and magnetic fields that are assumed to be spatially and temporally constant. 
Further, for the sake of simplicity, we consider a coordinate system where $q_{p}=m_{p}=c=1$.
Since all of our experiments have analytical solutions, the initial conditions of both the multi step polynomial and exponential methods are appropriately populated.
Finally, in each numerical example, the exponential PC method used is the 'PC1' stencil described in \cite{Glaser2009}. 
Specifically, we effected the SVD required to obtain the coefficients to a tolerance $10^{-12}$ (this serves as a lower bound for the error) as well as chose a semidisk radius $\rho=3.15$ and constructed a method that uses 22 steps of prior timestep information interpolated with 18 exponential basis sets.

\subsubsection{Linear Acceleration $\mathbf{E}(\mathbf{r},t)=(1,0,0)$, $\mathbf{B}(\mathbf{r},t)=\bar{0}$}
\begin{figure}
    \centering
    \includegraphics[width=0.5\textwidth]{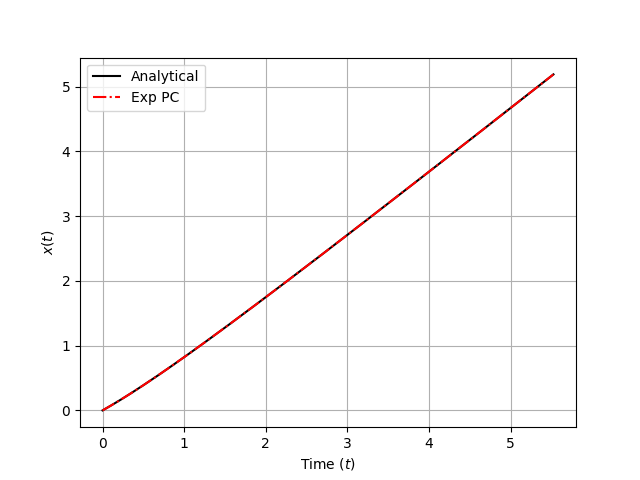}
    \caption{Comparison of particle trajectory between exponential PC and an analytical solution.}
    \label{fig:linear}
\end{figure}
\begin{figure}
    \centering
    \includegraphics[width=0.5\textwidth]{"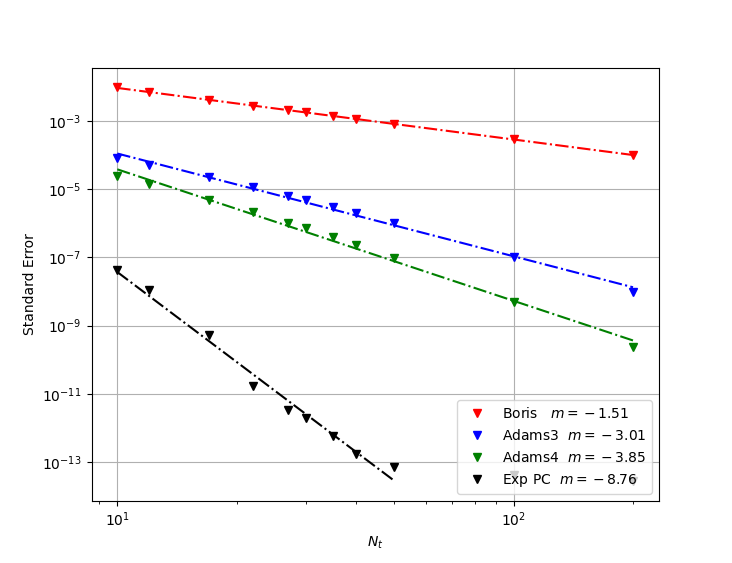"}
    \caption{Convergence of error (defined as the relative norm between computed and analytical particle trajectories) for the linear acceleration example. Note that the curve for the exponential PC is truncated due to hitting the floor of the tolerance of the method.}
    \label{fig:linear_conv}
\end{figure}
\begin{table}
    \centering
    \begin{tabular}{|c|c|c|}
        \hline
        Quantity & Natural & MKS  \\
        \hline
        $q_{p}$ & 1 & $1.6\times10^{-19}\ \text{C}$ \\
        $m_{p}$ & 1 & $9.1\times10^{-31}\ \text{kg}$ \\
        $c$ & 1 & $3\times 10^{8}\ \text{m/s}$ \\
        $\mathbf{E}_{p}(\mathbf{r},t)$ & 1 & $5.6\times 10^{-12}\ \text{V/m}$ \\
        $v_{0}$ & $1/\sqrt{2}$ & $c/\sqrt{2}\ \text{m/s}$ \\
        \hline
    \end{tabular}
    \caption{Conversion between natural and MKS units for the direct $\mathbf{E}$ acceleration example}
    \label{tab:unit_conversion_direct}
\end{table}
To begin, we consider the simple example of a particle being accelerated by an electric field of peak amplitude $1$ oriented along the $\hat{x}$ axis.
A mapping from natural to MKS units are provided in \ref{tab:unit_conversion_direct}.
The initial relativistic momentum of the particle was set such that $\mathbf{p}^{0}=\gamma_{p}^{0} \mathbf{v}_{p}^{0} = (1,0,0)$.
Fig. \ref{fig:linear}. shows a plot of the trajectories that the particles were predicted to follow by the exponential method, laid over the analytical result.
All of the methods closely agreed with the predicted curve, though the error in the position was significantly lower in the exponential method.
This is demonstrated in Fig. \ref{fig:linear_conv}, where we see Boris and the Adams predictor corrector methods roughly follow the trends expected from the polynomial order. The coarsest and finest timestep size used are $\Delta_{t}=0.1$ and $\Delta_{t}=0.005$ respectively. 
The exponential method, on the other hand, performs much better, showing better than 8th order convergence before flattening out due to hitting the tolerance of the method.
The analytical expression used is as follows, with the full derivation available at \cite{Landau1980Classical}:
\begin{equation}
    \mathbf{r}(t) = \frac{m}{q|E|}\left(\sqrt{1+\left(1+\frac{q|E|t}{m}\right)^2} - \gamma(0) \right).
\end{equation}
The factor $|E|$ refers to the magnitude of the electric field used, and $\gamma(0)$ refers to the initial relativistic factor.

\subsubsection{Cyclotron motion $\mathbf{E}(\mathbf{r},t)=\bar{0}$, $\mathbf{B}(\mathbf{r},t)=(0,0,1)$}
\begin{figure}
    \centering
    \includegraphics[width=0.5\textwidth]{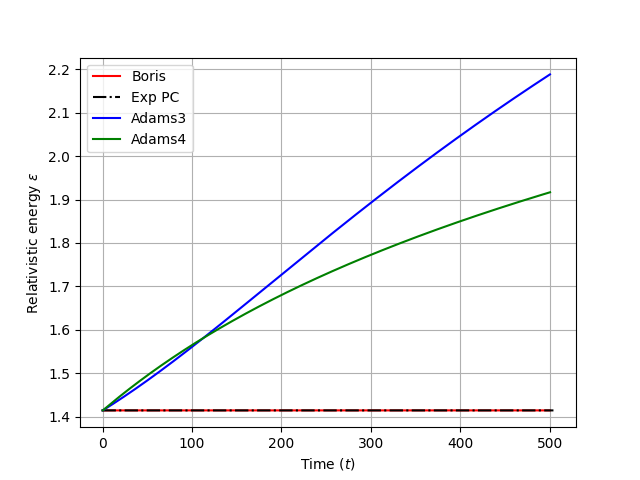}
    \caption{Comparison of particle energy over multiple cyclotron cycles. We note that the polynomial PC schemes spuriously gain energy for the timestep size chosen.}
    \label{fig:many_cycle_energy}
\end{figure}
\begin{figure}
    \centering
    \includegraphics[width=0.5\textwidth]{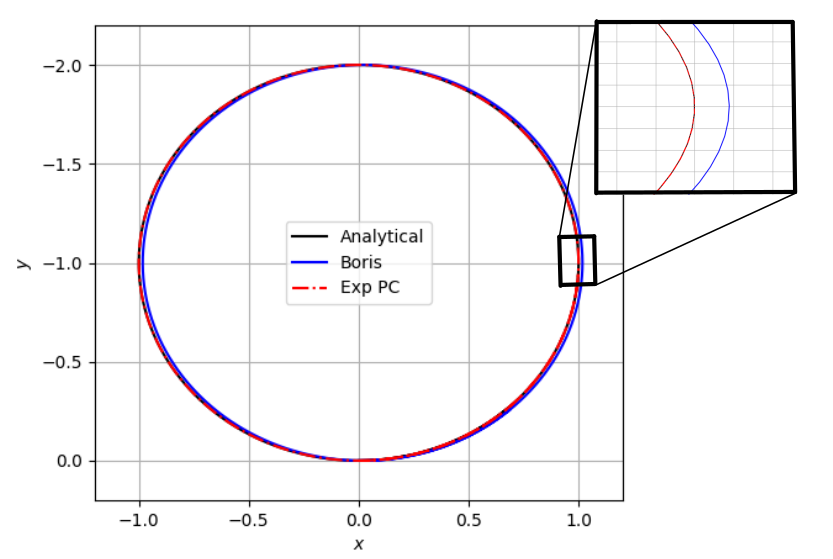}
    \caption{Comparison of particle trajectories obtained from the exponential method compared to Boris and the analytical solution. Note how the improved accuracy of the exponential method causes the trajectory to lie exactly on the analytical curve, whereas the Boris solution is slightly shifted for the same timestep size.}
    \label{fig:single_cycle_traj}
\end{figure}
\begin{figure}
    \centering
    \includegraphics[width=0.5\textwidth]{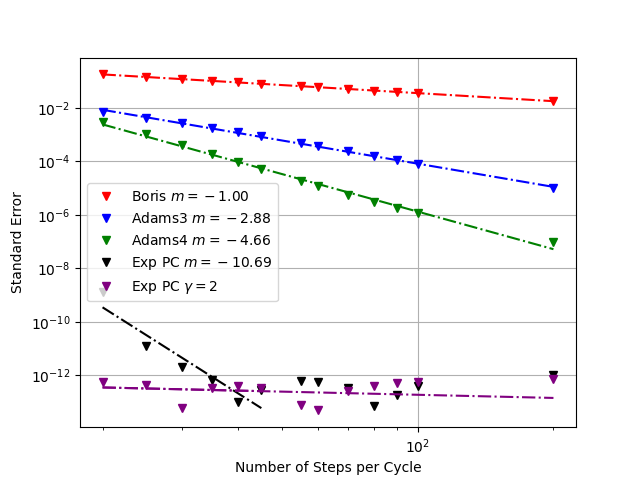}
    \caption{Convergence of position error (defined as the relative norm between computed and analytical particle trajectories)} over a single cycle for different evolution schemes.
    Once again, we note that the polynomial methods show approximately the expected convergence for the order of the polynomial used, while the exponential method converges much faster (11th order) before it hits the set tolerance.
    \label{fig:convergence_single_cycle}
\end{figure}
\begin{table}
    \centering
    \begin{tabular}{|c|c|c|}
        \hline
        Quantity & Natural & MKS  \\
        \hline
        $q_{p}$ & 1 & $1.6\times10^{-19}\ \text{C}$ \\
        $m_{p}$ & 1 & $9.1\times10^{-31}\ \text{kg}$ \\
        $c$ & 1 & $3\times 10^{8}\ \text{m/s}$ \\
        $\mathbf{B}_{p}(\mathbf{r},t)$ & 1 & $1.2\times 10^{-3}\ \text{T}$ \\
        $r_{c}$ & $1$ & $1\ \text{m}$ \\
        $v_{0}$ & $1/\sqrt{2}$ & $c/\sqrt{2}\ \text{m/s}$ \\
        \hline
    \end{tabular}
    \caption{Conversion between natural and MKS units for the cyclotron motion example}
    \label{tab:unit_conversion_cyclotron}
\end{table}
Next, we considered the dynamics of a particle moving under a constant $\hat{z}$-directed magnetic field.
As before, the initial momentum of the particle was set to $\mathbf{p}^{0}=\gamma_{p}^{0} \mathbf{v}_{p}^{0} = (1,0,0)$.
A mapping from natural to MKS units are provided in \ref{tab:unit_conversion_cyclotron}.
The canonical solution for this problem predicts a shortening of the non-relativistic cyclotron frequency and an increase in the radius of the loop by exactly a factor of $\gamma$.
The cyclotron test stresses a number of different aspects of the Newton solver.
First, since the particle is moving only under the influence of a magnetic field, this is the perfect setup to validate our predictions about energy conservation.
We note from Fig. \ref{fig:many_cycle_energy} that while Boris and the exponential method conserve the relativistic energy $m_{p}\gamma_{p}c$ to machine precision, the polynomial methods exhibit spurious heating.
Second, since the location of the gyrocenter is extremely sensitive to early time shifts, one can look at how well the different methods predict its location. As is evident from Fig. \ref{fig:single_cycle_traj}, despite conserving energy exactly, the Boris method has a slight shift of the orbit that is not seen in the more accurate exponential method.
We can also look at error convergence in the position over a single loop, as shown in Fig. \ref{fig:convergence_single_cycle}. The coarsest and finest stepsizes in the figure correspond to $\Delta_{t}=0.5$ and $\Delta_{t}=0.05$ respectively.
Here, we note that the Boris push and polynomial predictor corrector schemes converge at a slope approximately corresponding to the order of representation.
The exponential method on the other hand, shows much better convergence.
To stress the system further, the orbits were made more relativistic by setting the initial particle velocity such that $\gamma=2$. 
Looking once again at Fig. \ref{fig:convergence_single_cycle}, we see that the convergence is even better with the method hitting its tolerance with as little as 10 steps per cycle.
As in the previous case, the analytical solution is derived in \cite{Landau1980Classical}, and looks as follows:
\begin{equation}
    \mathbf{r}(t) = \frac{v_{0}\gamma_0}{|B|}\left[\sin\left(\frac{\gamma_{0}q|B|t}{m}\right)\hat{x} + \left(-1+\cos\left(\frac{\gamma_{0}q|B|t}{m}\right)\right)\hat{y}\right].
\end{equation}
The term $|B|$ refers to the magnitude of the $\hat{z}$ directed applied magnetic field.

\subsubsection{$|\mathbf{E}|=|\mathbf{B}|$}
\begin{figure}
    \centering
    \includegraphics[width=0.5\textwidth]{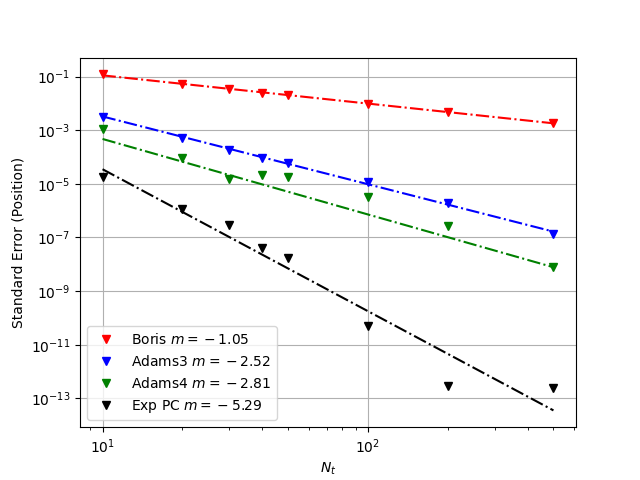}
    \caption{Convergence of error (defined as the relative norm between computed and analytical particle trajectories)} for the $|\mathbf{E}|=|\mathbf{B}|$ setup. The exponential method still outperforms Boris and the polynomial PC methods, but by a smaller order.
    \label{fig:convergence_guiding_center}
\end{figure}
\begin{table}
    \centering
    \begin{tabular}{|c|c|c|}
        \hline
        Quantity & Natural & MKS  \\
        \hline
        $q_{p}$ & 1 & $1.6\times10^{-19}\ \text{C}$ \\
        $m_{p}$ & 1 & $9.1\times10^{-31}\ \text{kg}$ \\
        $c$ & 1 & $3\times 10^{8}\ \text{m/s}$ \\
        $\mathbf{B}_{p}(\mathbf{r},t)$ & 1 & $1.2\times 10^{-3}\ \text{T}$ \\
        $\mathbf{E}_{p}(\mathbf{r},t)$ & 1 & $5.6\times 10^{-12}\ \text{V/m}$ \\
        $v_{0}$ & $1/\sqrt{2}$ & $c/\sqrt{2}\ \text{m/s}$ \\
        \hline
    \end{tabular}
    \caption{Conversion between natural and MKS units for the $|\mathbf{E}|=|\mathbf{B}|$ case. The conversion ensures that the electric and magnetic \emph{forces} acting on the positron at $t=0$ are equal.}
    \label{tab:unit_conversion_guiding}
\end{table}
Finally, we considered the more complex case of relativistic motion forced both by mutually perpendicular electric and magnetic fields of the same magnitude. 
As before, a mapping from natural to MKS units are provided in \ref{tab:unit_conversion_guiding}.
In this instance, we chose $\mathbf{E}(\mathbf{r},t)=(0,1,0)$ and $\mathbf{B}(\mathbf{r},t)=(0,0,1)$, with the particle starting at rest at the origin.
While the dynamics of the particle are more complicated as one reaches relativistic speeds, it is possible to derive an analytical solution, as done in \cite{Landau1980Classical}.
The final expression for our choice of parameters simplifies to
\begin{equation}
    \mathbf{r}(t) = \frac{U^{3}}{6}\hat{x} + \frac{U^{2}}{2}\hat{y},
\end{equation}
where
\begin{equation}
U = \frac{\left(\left(\sqrt{9t^2+8} + 3t\right)^{2/3} - 2\right)}{\left(\sqrt{9t^2+8}+3t\right)^{1/3}}.
\end{equation}
Upon running this setup with each of our solvers, we once again observe from Fig. \ref{fig:convergence_guiding_center} that the exponential method significantly outperforms the polynomial based methods, though by a smaller polynomial factor than the other examples presented thus far.
The coarsest and finest stepsizes in Fig. 6 correspond to $\Delta_{t}=1.0$ and $\Delta_{t}=0.02$ respectively.

\subsection{Space Varying Fields}
\begin{figure}
    \centering
    \includegraphics[width=0.5\textwidth]{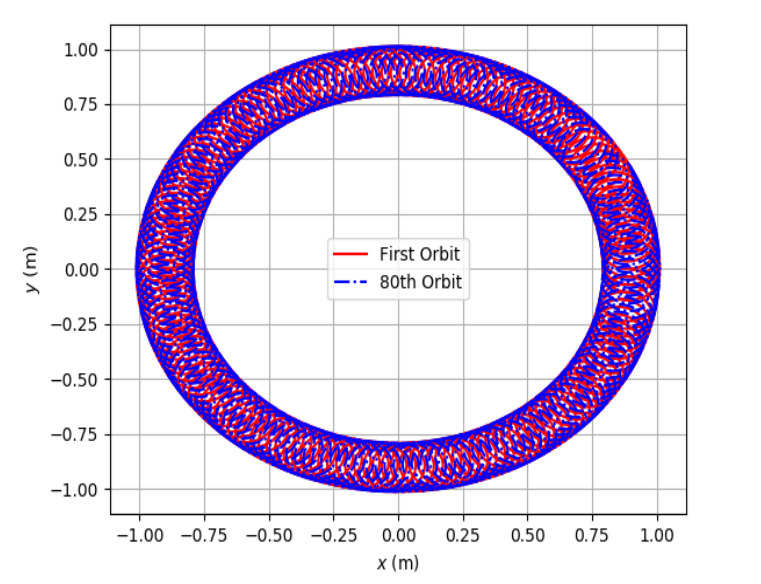}
    \caption{Plot of particle trajectory for the first and 80th orbits, showing good preservation of phase space volume}
    \label{fig:orbits}
\end{figure}
Next, we tested the long-time stability of the exponential method by considering a system wherein the scalar potential and magnetic fields varied as follows, $\phi(\mathbf{r},t) = 0.01\sqrt{x^{2} + y^{2}}$, $\mathbf{B}(\mathbf{r},t) = \hat{z}\sqrt{x^{2} + y^{2}}$. 
Once again employing units of $m_{p}=q_{p}=c=1$, the particle was initialized at $\mathbf{r}_{p}^{0} = (0.9,0,0)$ and $\mathbf{u}_{p}^{0}=(0.1,0,0)$.
Further, we chose a step size $\Delta_{t}=1.5\times 10^{-5}$ and ran the system for 50000 timesteps.
We note from Fig. \ref{fig:orbits} that exponential PC method accurately tracks both the fast gyrations and the slow rotations without losing or gaining energy.

\subsection{Particle trapped in a Magnetic Bottle}
\begin{figure}
    \centering
    \includegraphics[width=0.5\textwidth]{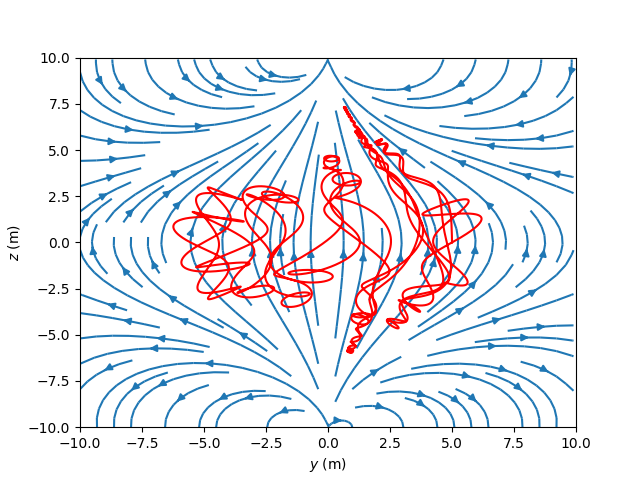}
    \caption{Plot of particle trajectories overlaid on a lineplot of the magnetic field along the $y$-$z$ plane.}
    \label{fig:bottle_fields}
\end{figure}
\begin{figure}
    \centering
    \includegraphics[width=0.5\textwidth]{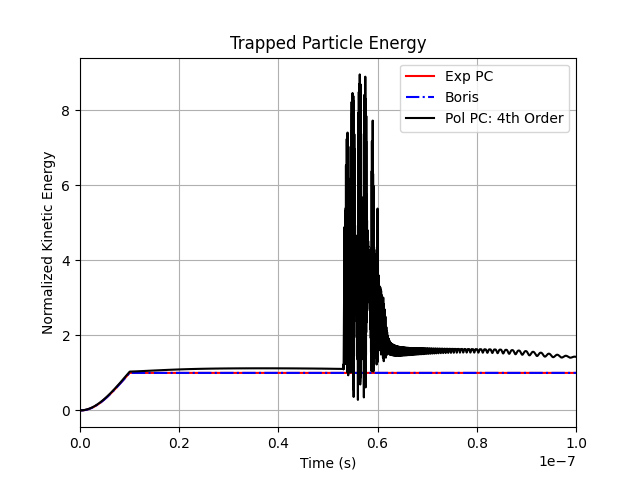}
    \caption{Comparison of particle energy for different particle update schemes. The Kinetic energy is normalized to the energy at the point where the electric field is switched off.}
    \label{fig:energy_bottle}
\end{figure}
To stress energy conservation further, we considered a setup where an electron was initialized within a magnetic bottle defined using two magnetic dipole oriented along the $z$-direction, with field strength
\begin{equation}
    \mathbf{B}(\mathbf{r}) = \frac{\mu_{0}}{4\pi}\left(\frac{3\mathbf{r}(\mathbf{m}\cdot\mathbf{r})}{|\mathbf{r}|^{5}} - \frac{\mathbf{m}}{|\mathbf{r}|^{3}}\right).
\end{equation}
The two dipoles were separated by a distance of 20 m and the magnetic moment was set to $m=10^{5}\hat{z}$.
The mass, charge and speed of light were represented in MKS units.
The particle initially begins at rest at position $(0,5,0)$ and is accelerated for 10 ns by a $y$-directed electric field before it is switched off and moves only in the magnetic field from that point.
The field lines of the dipole, as well as the trapped trajector predicted by the exponential method are shown in Fig. \ref{fig:bottle_fields}.
We note from Fig. \ref{fig:energy_bottle} that the exponential method and Boris once again conserve energy once the electric field is switched off.
The polynomial method, on the other hand, spuriously gains energy after the electric field is switched off.
It also generates massive spikes in energy when the particle is near the deflection point in its trajectory, due to the method not resolving large changes in $\mathbf{v}(t)\times \mathbf{B}(\mathbf{r},t)$ in a manner that preserves kinetic energy. 

\subsection{Expanding Particle Beam}
\begin{figure}
    \centering
    \includegraphics[width=0.5\textwidth]{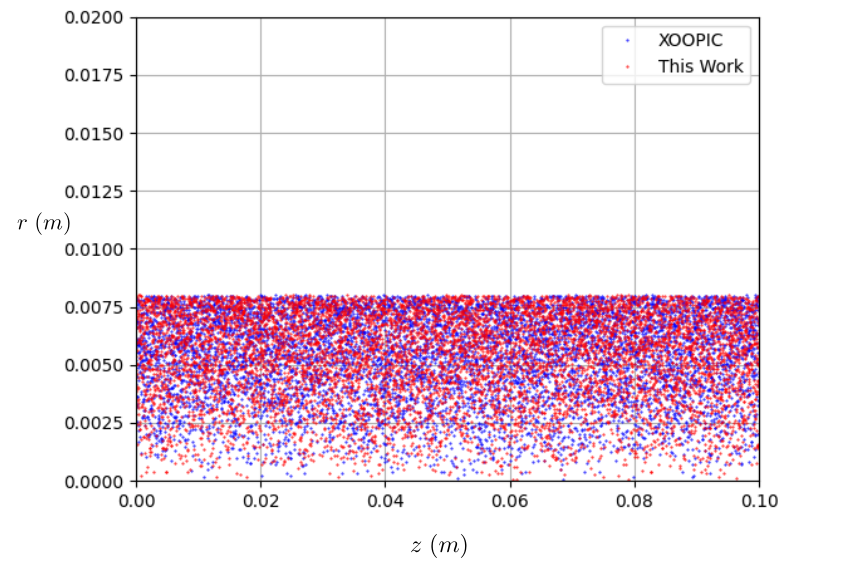}
    \caption{Snapshot of the particle spread at $t=10\ \text{ns}$. We note that
    the particle spead generally follows the expansion expected from the 
    physics, and further agrees well with equivalent results from XOOPIC}
    \label{fig:beam_snap}
\end{figure}
\begin{figure}
    \centering
    \includegraphics[width=0.5\textwidth]{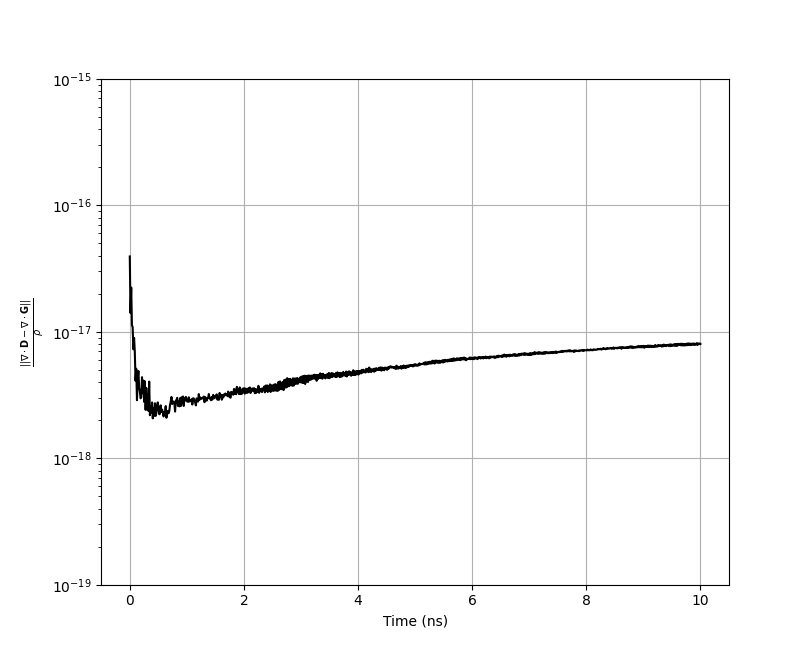}
    \caption{Satisfaction of Gauss' Law over the run for the 
    expanding particle beam. The figure shown is the relative norm between the $\nabla\cdot \mathbf{D}(\mathbf{r},t)$ and $\nabla\cdot \mathbf{G}(\mathbf{r},t)$.}
    \label{fig:beam_gauss}
\end{figure}
\begin{figure}
    \centering
    \includegraphics[width=0.5\textwidth]{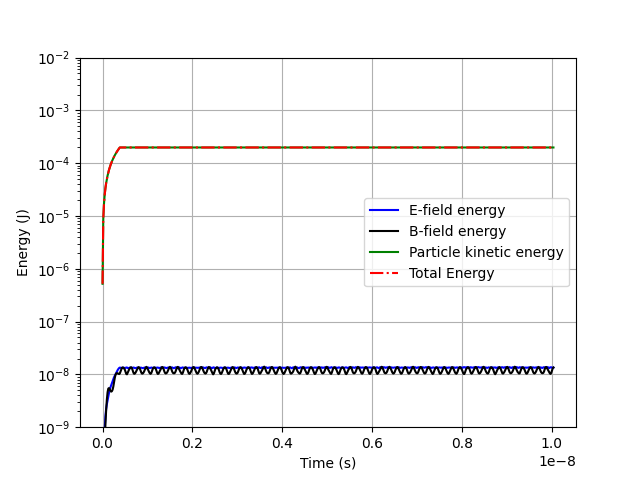}
    \caption{Time history of particle and field energy for the expanding beam. Due to the high initial particle momentum, the energy is largely dominated by the kinetic energy of the particles. However, we note from the data that both field and particle energy are conserved across the run.}
    \label{fig:beam_energy}
\end{figure}
\begin{figure}
    \centering
    \includegraphics[width=0.5\textwidth]{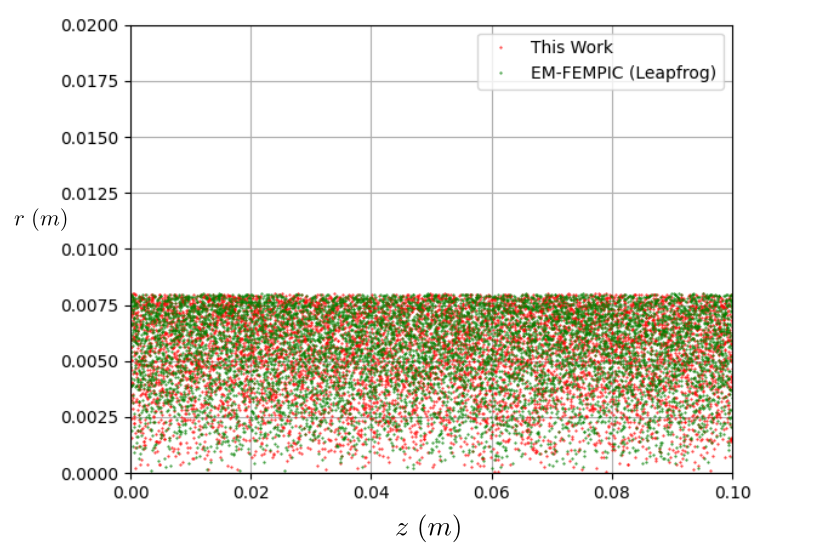}
    \caption{Snapshot of the particle spread at $t=10\ \text{ns}$ predicted by the proposed method compared to results from an EM-FEMPIC formulation run using a leapfrog integrator. Again, we note very good agreement between the two methods.}
    \label{fig:beam_snap_lf}
\end{figure}
To test the validity of the exponential predictor-corrector scheme on a full
EM-FEMPIC scheme, we analysed the physics of an expanding beam contained within
a conducting cylindrical cavity of length 10 cm and radius 2 cm oriented such 
that the axis of rotation aligned along $\hat{z}$.
The walls of the cavity were assumed to be perfectly conducting and a particle
beam composed of electrons was initialized at $z=0$.
The beam was characterized by a driving voltage of $500\ kV$, a driving 
current of $1\ A$ and an emitter radius of $0.8$ cm.
This yielded an initial setup wherein the electrons were accelerated at roughly $\gamma=2$ or $\mathbf{v}_{p}^{0}\sim 2.6\times 10^{8}$ m/s.
The exponential predictor-corrector was used with a timestep size equal to $\Delta_{t}=15$ ps -- while the comparison runs with XOOPIC\cite{xoopic} and an explicit EM-FEMPIC code were run at $\Delta_{t}=1$ ps.
The tetrahedral mesh used to discretize the system had an average cell size of 8.2 mm and the maximum number of particles within the system at any one time was approximately 8000 particles, resulting in approximately 3.5 particles per cell.
At each step, one predictor step was taken followed by as many corrector steps required as to reach the desired tolerance of $10^{-9}$.
In this specific run, achieving this threshold rarely required more than one corrector step.
The particle beam also had a turn-on time of 
1 ns so as to have a smooth increase in field strength.
A snapshot of of the particle beam at $10$ ns -- compared with an equivalent simulation set-up on XOOPIC\cite{xoopic} -- is shown in Fig. 
\ref{fig:beam_snap} and the satisfaction of Gauss' Law over the course of the
entire run is shown in Fig. \ref{fig:beam_gauss}.
As is evident, the proposed method both agrees with the expected physics and
satisfies Gauss' Laws to machine precision.
Likewise, we note from Fig. \ref{fig:beam_energy} that particle and field energies are conserved over the course of the simulation.
Finally, we note from Fig. \ref{fig:beam_snap_lf} that the proposed method compares well to an explicit leapfrog based EM-FEMPIC solver, albeit at a stepsize that was 15 times larger.

\subsection{Particle Beam in a Klystron}
\begin{figure}
    \centering
    \includegraphics[width=0.5\textwidth]{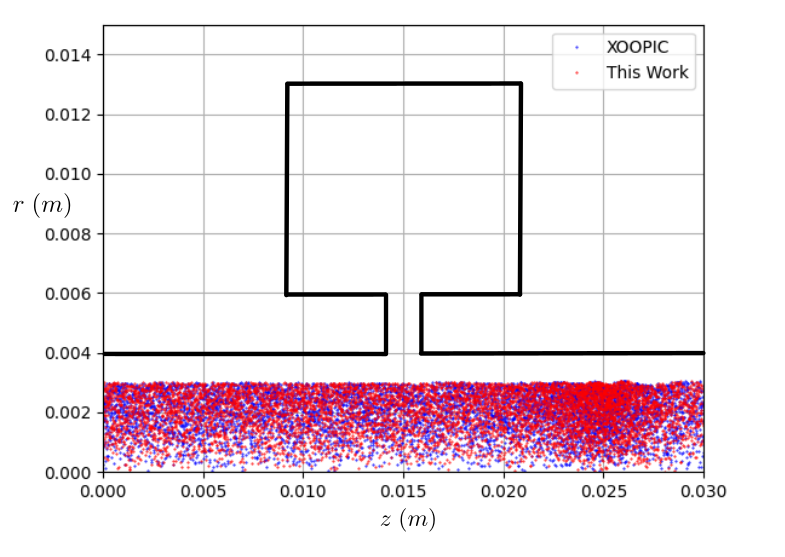}
    \caption{Snapshot of the position distribution at $t=4\ \text{ns}$. The proposed method agrees well with XOOPIC and predicts the expected particle bunching behavior around the same region as XOOPIC (at approximately $z=2.5$ cm).}
    \label{fig:klystron_dist}
\end{figure}
\begin{figure}
    \centering
    \includegraphics[width=0.5\textwidth]{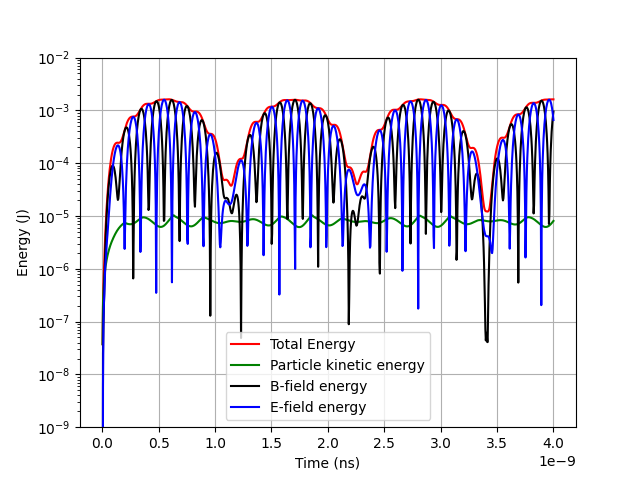}
    \caption{Time history of particle and field energy for the klystron. Unlike in the previous run, most of the energy is in the fields -- we note here that once again average energy is conserved across time for both fields and particles.}
    \label{fig:klystron_energy}
\end{figure}
To stress the proposed solver in a system with excited cavity modes, we analysed the behavior of a simple cylindrical klystron, with dimensions as in Fig. \ref{fig:klystron_dist}. Once again, the system was discretized by a tetrahedral mesh with average edge length of 2.93 mm. The solver was run with a step size of 10 ps and compared to XOOPIC running at $1$ ps. 
A sinusoidal current source was defined on a cylindrical surface existing at the middle of the neck of the klystron, with an amplitude of 100 A and frequency of $4.1$ GHz. These dimensions correspond to an average of around 2 particles per cell, though there are regions in the domain with higher density as shown by the bunching observed in Fig. \ref{fig:klystron_dist}. Once again, we note good agreement with XOOPIC, and observe from Fig. \ref{fig:klystron_energy} that both field and particle energies are conserved on average over multiple cycles of the run.

\section{Summary \label{sec:summary}}
In this work, we have proposed and validated a fully relativistic charge 
conserving EM-FEMPIC scheme that makes use of a novel exponential predictor
corrector framework. 
We have numerically benchmarked our proposed framework over both single particle full PIC runs, demonstrating superior accuracy in the former and good agreement with existing validated methods in the latter.
The improvement of scalability of this scheme to make it more applicable to 
complex devices that are of interest in high power design is currently being 
analysed and will be reported on in the near future.

\section*{Acknowledgments}
The authors would like to thank the HPCC Facility, Michigan State University, East Lansing, MI, USA. The authors would also like to acknowledge numerous discussions with Dr. Zane D. Crawford on charge conserving EM-FEMPIC.

\section{Appendix}

\typeout{}
\bibliographystyle{IEEEtran}
\bibliography{HOtime}

\begin{thebibliography}{10}
\providecommand{\url}[1]{#1}
\csname url@samestyle\endcsname
\providecommand{\newblock}{\relax}
\providecommand{\bibinfo}[2]{#2}
\providecommand{\BIBentrySTDinterwordspacing}{\spaceskip=0pt\relax}
\providecommand{\BIBentryALTinterwordstretchfactor}{4}
\providecommand{\BIBentryALTinterwordspacing}{\spaceskip=\fontdimen2\font plus
\BIBentryALTinterwordstretchfactor\fontdimen3\font minus
  \fontdimen4\font\relax}
\providecommand{\BIBforeignlanguage}[2]{{%
\expandafter\ifx\csname l@#1\endcsname\relax
\typeout{** WARNING: IEEEtran.bst: No hyphenation pattern has been}%
\typeout{** loaded for the language `#1'. Using the pattern for}%
\typeout{** the default language instead.}%
\else
\language=\csname l@#1\endcsname
\fi
#2}}
\providecommand{\BIBdecl}{\relax}
\BIBdecl

\bibitem{marchand2011ptetra}
R.~Marchand, ``{PTETRA}, a tool to simulate low orbit satellite--plasma
  interaction,'' \emph{IEEE Transactions on Plasma Science}, vol.~40, no.~2,
  pp. 217--229, 2011.

\bibitem{lemke1999three}
R.~Lemke, T.~Genoni, and T.~Spencer, ``Three-dimensional particle-in-cell
  simulation study of a relativistic magnetron,'' \emph{Physics of Plasmas},
  vol.~6, no.~2, pp. 603--613, 1999.

\bibitem{fourkal2002particle}
E.~Fourkal, B.~Shahine, M.~Ding, J.~Li, T.~Tajima, and C.-M. Ma, ``Particle in
  cell simulation of laser-accelerated proton beams for radiation therapy,''
  \emph{Medical Physics}, vol.~29, no.~12, pp. 2788--2798, 2002.

\bibitem{pinto2014charge}
M.~C. Pinto, S.~Jund, S.~Salmon, and E.~Sonnendr{\"u}cker, ``Charge-conserving
  fem--pic schemes on general grids,'' \emph{Comptes Rendus Mecanique}, vol.
  342, no. 10-11, pp. 570--582, 2014.

\bibitem{moon2014exact}
H.~Moon, F.~L. Teixeira, and Y.~A. Omelchenko, ``Exact charge-conserving
  scatter--gather algorithm for particle-in-cell simulations on unstructured
  grids: A geometric perspective,'' \emph{Computer Physics Communications},
  vol. 194, pp. 43--53, 2015.

\bibitem{squire2012geometric}
J.~Squire, H.~Qin, and W.~M. Tang, ``Geometric integration of the
  vlasov-maxwell system with a variational particle-in-cell scheme,''
  \emph{Physics of Plasmas}, vol.~19, no.~8, p. 084501, 2012.

\bibitem{moon2015exact}
H.~Moon, F.~L. Teixeira, and Y.~A. Omelchenko, ``Exact charge-conserving
  scatter--gather algorithm for particle-in-cell simulations on unstructured
  grids: A geometric perspective,'' \emph{Computer Physics Communications},
  vol. 194, pp. 43--53, 2015.

\bibitem{o2021set}
S.~O'Connor, Z.~Crawford, J.~Verboncoeur, J.~Lugisland, and B.~Shanker, ``A set
  of benchmark tests for validation of 3d particle in cell methods,''
  \emph{arXiv preprint arXiv:2101.09299}, 2021.

\bibitem{tonti}
E.~Tonti, ``Finite formulation of electromagnetic field,'' \emph{IEEE
  Transactions on Magnetics}, vol.~38, no.~2, pp. 333--336, 2002.

\bibitem{deschamps}
G.~Deschamps, ``Electromagnetics and differential forms,'' \emph{Proceedings of
  the IEEE}, vol.~69, no.~6, pp. 676--696, 1981.

\bibitem{kormann2021energy}
K.~Kormann and E.~Sonnendr{\"u}cker, ``Energy-conserving time propagation for a
  structure-preserving particle-in-cell vlasov--maxwell solver,'' \emph{Journal
  of Computational Physics}, vol. 425, p. 109890, 2021.

\bibitem{xiao_2018}
\BIBentryALTinterwordspacing
J.~XIAO, H.~QIN, and J.~LIU, ``Structure-preserving geometric particle-in-cell
  methods for vlasov-maxwell systems,'' \emph{Plasma Science and Technology},
  vol.~20, no.~11, p. 110501, sep 2018. [Online]. Available:
  \url{https://doi.org/10.1088/2058-6272/aac3d1}
\BIBentrySTDinterwordspacing

\bibitem{CamposPinto2022}
\BIBentryALTinterwordspacing
M.~Campos~Pinto, K.~Kormann, and E.~Sonnendr{\"u}cker, ``Variational framework
  for structure-preserving electromagnetic particle-in-cell methods,''
  \emph{Journal of Scientific Computing}, vol.~91, no.~2, p.~46, Mar 2022.
  [Online]. Available: \url{https://doi.org/10.1007/s10915-022-01781-3}
\BIBentrySTDinterwordspacing

\bibitem{crawford2021rubrics}
Z.~D. Crawford, S.~O'Connor, J.~Luginsland, and B.~Shanker, ``Rubrics for
  charge conserving current mapping in finite element particle in cell
  methods,'' \emph{arXiv preprint arXiv:2101.12128}, 2021.

\bibitem{Oconnor_time_integration}
S.~O'Connor, Z.~D. Crawford, O.~H. Ramachandran, J.~Luginsland, and B.~Shanker,
  ``Time integrator agnostic charge conserving finite element pic,''
  \emph{Physics of Plasmas}, vol.~28, no.~9, p. 092111, 2021.

\bibitem{oconnor2021quasihelmholtz}
S.~O'Connor, Z.~D. Crawford, O.~Ramachandran, J.~Luginsland, and B.~Shanker,
  ``Quasi helmholtz decomposition, gauss laws and charge conservation for
  finite element pic,'' \emph{Computer Physics Communications}, vol. 276, p.
  108345, 2022.

\bibitem{ramachandran2022envelopetracking}
\BIBentryALTinterwordspacing
O.~H. Ramachandran, Z.~D. Crawford, S.~O'Connor, J.~Luginsland, and B.~Shanker,
  ``An envelope tracking approach for particle in cell simulations,'' 2022.
  [Online]. Available: \url{https://arxiv.org/abs/2208.12795}
\BIBentrySTDinterwordspacing

\bibitem{boris}
J.~P. Boris, ``Relativistic plasma simulation-optimization of a hybrid code,''
  \emph{Proceeding of Fourth Conference on Numerical Simulations of Plasmas},
  November 1970.

\bibitem{why_boris_is_good}
\BIBentryALTinterwordspacing
H.~Qin, S.~Zhang, J.~Xiao, J.~Liu, Y.~Sun, and W.~M. Tang, ``Why is boris
  algorithm so good?'' \emph{Physics of Plasmas}, vol.~20, no.~8, p. 084503,
  2013. [Online]. Available: \url{https://doi.org/10.1063/1.4818428}
\BIBentrySTDinterwordspacing

\bibitem{2018_boris_alt_method}
\BIBentryALTinterwordspacing
S.~Zenitani and T.~Umeda, ``On the boris solver in particle-in-cell
  simulation,'' \emph{Physics of Plasmas}, vol.~25, no.~11, p. 112110, 2018.
  [Online]. Available: \url{https://doi.org/10.1063/1.5051077}
\BIBentrySTDinterwordspacing

\bibitem{Bacchini_2019}
\BIBentryALTinterwordspacing
F.~Bacchini, J.~Amaya, and G.~Lapenta, ``The relativistic implicit
  particle-in-cell method,'' \emph{Journal of Physics: Conference Series}, vol.
  1225, no.~1, p. 012011, may 2019. [Online]. Available:
  \url{https://doi.org/10.1088/1742-6596/1225/1/012011}
\BIBentrySTDinterwordspacing

\bibitem{markidis_2011}
\BIBentryALTinterwordspacing
S.~Markidis and G.~Lapenta, ``The energy conserving particle-in-cell method,''
  \emph{Journal of Computational Physics}, vol. 230, no.~18, pp. 7037--7052,
  2011. [Online]. Available:
  \url{https://www.sciencedirect.com/science/article/pii/S0021999111003445}
\BIBentrySTDinterwordspacing

\bibitem{implicit_thesis}
S.~Markidis, \emph{Development of Implicit Kinetic Simulation Methods, and
  their Application to Ion Beam Propagation in Current and Future Neutralized
  Drift Compression Experiments}.\hskip 1em plus 0.5em minus 0.4em\relax
  University of Illinois at Urbana-Champaign ProQuest Dissertations Publishing,
  2010.

\bibitem{Glaser2009}
\BIBentryALTinterwordspacing
A.~Glaser and V.~Rokhlin, ``A new class of highly accurate solvers for ordinary
  differential equations,'' \emph{Journal of Scientific Computing}, vol.~38,
  no.~3, pp. 368--399, Mar 2009. [Online]. Available:
  \url{https://doi.org/10.1007/s10915-008-9245-1}
\BIBentrySTDinterwordspacing

\bibitem{zienkiewicz1977new}
O.~C. Zienkiewicz, ``A new look at the newmark, houbolt and other time stepping
  formulas. a weighted residual approach,'' \emph{Earthquake Engineering \&
  Structural Dynamics}, vol.~5, no.~4, pp. 413--418, 1977.

\bibitem{Boss88}
A.~Bossavit, ``Whitney forms: a class of finite elements for three-dimensional
  computations in electromagnetism,'' \emph{IEE Proceedings A - Physical
  Science, Measurement and Instrumentation, Management and Education -
  Reviews}, vol. 135, no.~8, pp. 493--500, Nov 1988.

\bibitem{crawford2020unconditionally}
Z.~Crawford, J.~Li, A.~Christlieb, and B.~Shanker, ``Unconditionally stable
  time stepping method for mixed finite element maxwell solvers,''
  \emph{Progress In Electromagnetics Research}, vol. 103, pp. 17--30, 2020.

\bibitem{Wong95}
M.-F. Wong, O.~Picon, and V.~F. Hanna, ``A finite element method based on
  whitney forms to solve maxwell equations in the time domain,'' \emph{IEEE
  Transactions on Magnetics}, vol.~31, no.~3, pp. 1618--1621, May 1995.

\bibitem{He06}
B.~He and F.~Teixeira, ``Geometric finite element discretization of maxwell
  equations in primal and dual spaces,'' \emph{Physics Letters A}, vol. 349,
  no. 1–4, pp. 1 -- 14, 2006.

\bibitem{Kett99}
A.~Bossavit and L.~Kettunen, ``Yee-like schemes on a tetrahedral mesh, with
  diagonal lumping,'' \emph{International Journal of Numerical Modelling:
  Electronic Networks, Devices and Fields}, vol.~12, no. 1-2, pp. 129--142,
  1999.

\bibitem{wang2010application}
R.~Wang, D.~J. Riley, and J.-M. Jin, ``Application of tree-cotree splitting to
  the time-domain finite-element analysis of electromagnetic problems,''
  \emph{IEEE transactions on antennas and propagation}, vol.~58, no.~5, pp.
  1590--1600, 2010.

\bibitem{venkatarayalu2006removal}
N.~V. Venkatarayalu and J.-F. Lee, ``Removal of spurious dc modes in edge
  element solutions for modeling three-dimensional resonators,'' \emph{IEEE
  transactions on microwave theory and techniques}, vol.~54, no.~7, pp.
  3019--3025, 2006.

\bibitem{Landau1980Classical}
L.~D. Landau and E.~M. Lifshitz, \emph{The Classical Theory of Fields}.\hskip
  1em plus 0.5em minus 0.4em\relax Butterworth-Heinemann, Jan. 1980.

\bibitem{xoopic}
\BIBentryALTinterwordspacing
J.~Verboncoeur, A.~Langdon, and N.~Gladd, ``An object-oriented electromagnetic
  pic code,'' \emph{Computer Physics Communications}, vol.~87, no.~1, pp.
  199--211, 1995, particle Simulation Methods. [Online]. Available:
  \url{https://www.sciencedirect.com/science/article/pii/001046559400173Y}
\BIBentrySTDinterwordspacing

\end{thebibliography}

\begin{IEEEbiography}[{\includegraphics[width=1in,height=1.25in,clip,keepaspectratio]{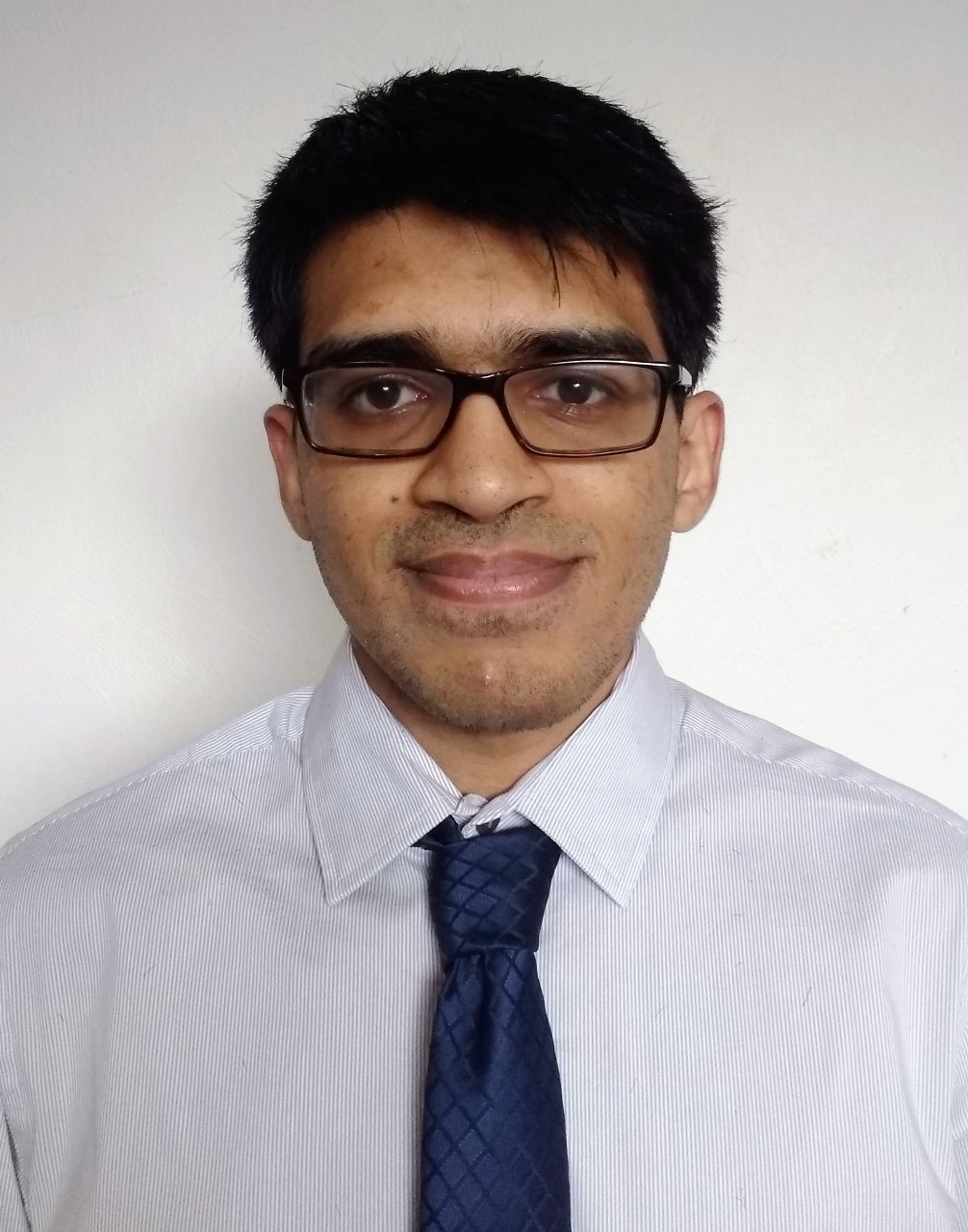}}]{Omkar H. Ramachandran} (S'20) received the B.A. degree
in physics from the University of Colorado at Boulder in 2018 and is currently pursuing the Ph.D in electrical and computer engineering at Michigan State University, East Lansing, MI. His research interests include several topics in computational electromagnetics, including particle-in-cell methods, analysis of coupled EM-device systems and nonlinear optimization.
\end{IEEEbiography}

\begin{IEEEbiography}[{\includegraphics[width=1in,height=1.25in,clip,keepaspectratio]{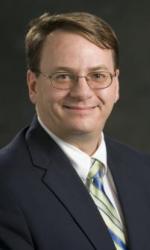}}]{Leo Kempel}
(S’89–M’94–SM’99–F’09) was
born in Akron, OH, USA, in 1965. He received the
B.S.E.E. degree from the University of Cincinnati,
Cincinnati, OH, USA, in 1989, and the M.S.E.E.
and Ph.D. degrees from the University of Michigan, Ann Arbor, MI, USA, in 1990 and 1994,
respectively.
After a brief post-doctoral appointment at the
University of Michigan, he joined Mission Research
Corporation, Goleta, CA, USA, in 1994, as a Senior
Research Engineer. He led several projects involving
the design of conformal antennas, computational electromagnetics, scattering
analysis, and high-power/ultrawideband microwaves. He joined Michigan
State University, East Lansing, MI, USA, in 1998. He served as an IPA
with the Air Force Research Laboratory’s Sensors Directorate, Riverside,
OH, USA, from 2004 to 2005 and 2006 to 2008. He was the Inaugural
Director of the Michigan State University High Performance Computing
Center, East Lansing. He was the first Associate Dean for Special Initiatives
with the College of Engineering, Michigan State University, from 2006 to
2008, and the Associate Dean for Research from 2008 to 2013. He then
became the Acting Dean of Engineering in 2013. Since 2014, he has been
the Dean with the College of Engineering, Michigan State University. He has
co-authored the book The Finite Element Method for Electromagnetics (IEEE
Press). His current research interests include computational electromagnetics,
conformal antennas, microwave/millimeter-wave materials, and measurement
techniques.
Prof. Kempel is a fellow of the Applied Computational Electromagnetics
Society (ACES). He was a member of the Antennas and Propagation Society’s
Administrative Committee and the ACES Board of Directors. He is a member
of Tau Beta Pi, Eta Kappa Nu, and Commission B of URSI. He served as
the Technical Chairperson for the 2001 ACES Conference and the Technical
Co-Chair for the Finite Element Workshop held in Chios, Greece, in 2002.
He was the Fellow Evaluation Committee Chairperson for the IEEE Antennas
and Propagation Society and served on the IEEE Fellow Board from 2013 to
2015. He was a recipient of the CAREER Award by the National Science
Foundation, the Teacher-Scholar Award by Michigan State University in
2002, and the MSU College of Engineering’s Withrow Distinguished Scholar
(Junior Faculty) Award in 2001. He served on the U.S. Air Force Scientific
Advisory Board from 2011 to 2015. He served as an Associate Editor of
the IEEE TRANSACTIONS ON ANTENNAS AND PROPAGATION. He is an
active reviewer for several IEEE publications as well as the Journal of
Electromagnetic Waves and Applications and Radio Science.
\end{IEEEbiography}

\begin{IEEEbiography}[{\includegraphics[width=1in,height=1.25in,clip,keepaspectratio]{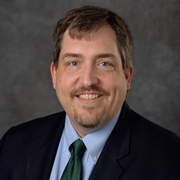}}]{John Luginsland} (Fellow, IEEE) received the B.S.E., M.S.E., and Ph.D. degrees in nuclear engineering from the University of Michigan, Ann Arbor, MI, USA, in 1992, 1994, and 1996, respectively. He served as a Professor with the Department of Computational Mathematics, Science, and Engineering and Electrical and Computer Engineering, Michigan State University, East Lansing, MI, USA, and various roles with the Air Force Office of Scientific Research, Arlington, VA, USA, including Acting Division Chief, Division Technical Advisor, Acting Branch Chief, Program Manager for Plasma Physics, and Program Manager for Laser Science. He is currently a Program Officer at AFOSR where he has previously served as the Program Element Monitor for Air Force Basic Research in the office of the Assistant Secretary of the Air Force for Acquisition. He was a Staff Member with NumerEx LLC, Ithaca, NY, USA, Science Applications International Corporation, Albuquerque, NM, USA, and the Air Force Research Laboratory, Kirtland AFB, NM, where he was also a National Research Council Post-Doctoral Researcher. He was a Senior Scientist and Principal Investigator with Confluent Sciences, LLC, Albuquerque. He was also an Adjunct Professor of electrical and computer engineering with Michigan State University, and a member of the Intelligence Science and Technology Experts Group (ISTEG) of the National Academy of Sciences, Washington, DC, USA, and works with the Board on Army Research and Development of the National Academy of Sciences. His research interests are in accelerator design, coherent radiation sources, dense kinetic plasmas, laser physics, serious games including agent-based models and wargames, as well as computational modeling using high-performance computing and machine learning techniques. Dr. Luginsland is a fellow of the Air Force Research Laboratory. He received the IEEE Nuclear and Plasma Science Society’s Early Achievement Award. He is currently a Vice-Chair of the IEEE Plasma Science and Applications Committee, a Chair of the IEEE PSAC, and a Guest Editor of the IEEE TRANSACTIONS ON PLASMA SCIENCE SPECIAL ISSUE ON HIGH POWER MICROWAVE SOURCES .
\end{IEEEbiography}

\begin{IEEEbiography}[{\includegraphics[width=1in,height=1.25in,clip,keepaspectratio]{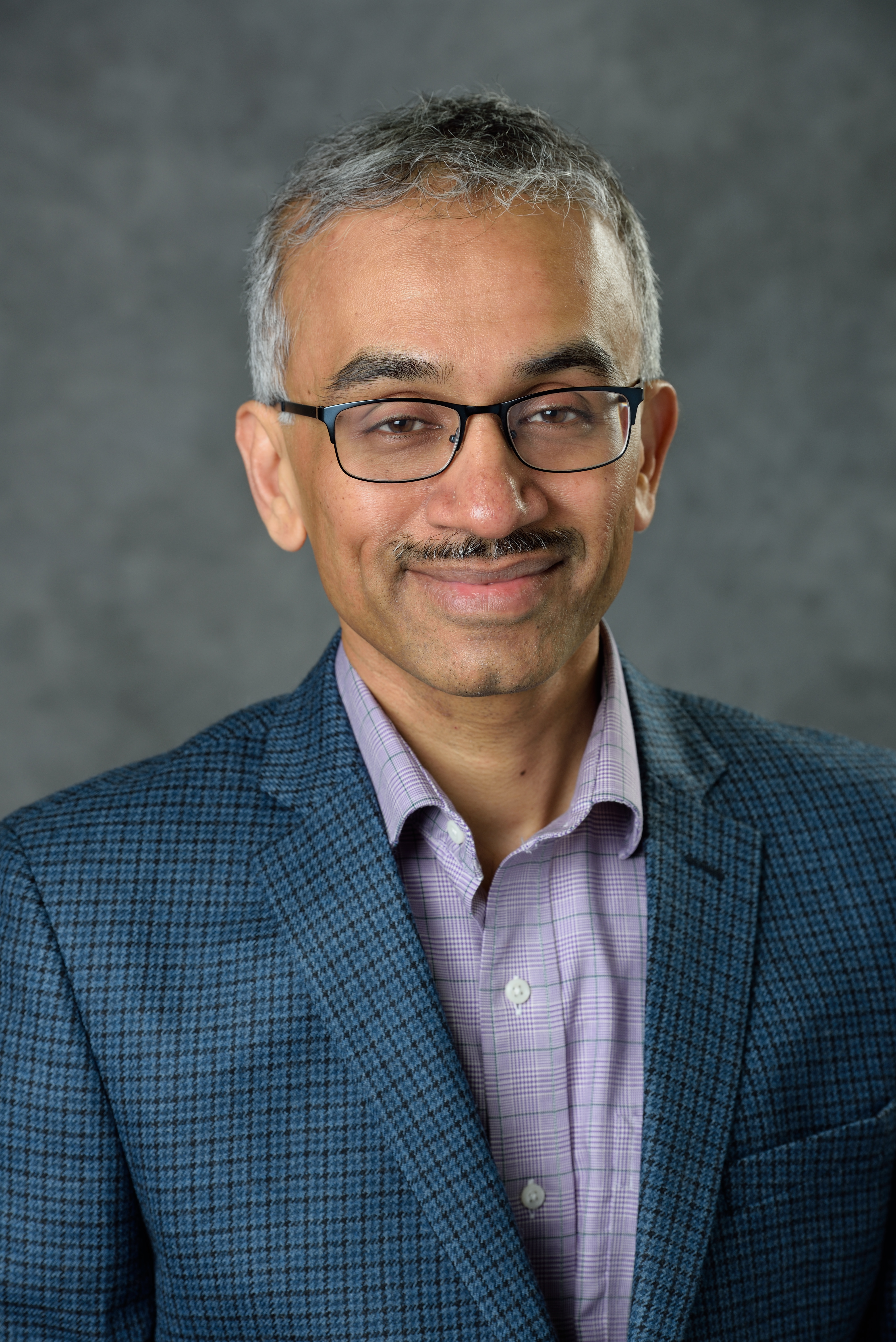}}]{B. Shanker} (Fellow, IEEE) received the B'Tech from the Indian Institute of Technology, Madras, India in 1989, M.S. and Ph.D in 1992 and 1993, respectively, from The Pennsylvania State University. From 1993 to 1996 he was a research associate in the Department of Biochemistry and Biophysics at Iowa State University where he worked on the Molecular Theory of Optical Activity. From 1996 to 1999 he was with the Center for Computational Electromagnetics at the University of Illinois at Urbana-Champaign as a Visiting Assistant Professor, and from 1999-2002 with the Department of Electrical and Computer Engineering at Iowa State University as an Assistant Professor. From 2017, he was a University Distinguished Professor (an honor accorded to about 2\% of tenure system MSU faculty members) in the Department of Electrical and Computer Engineering at Michigan State University, and the Department of Physics and Astronomy. Currently, he is a Professor and Chair of Electrical and Computer Engineering at The Ohio State University. At Michigan State University, he was appointed Associate Chair of the Department of Computational Mathematics, Science and Engineering, a new department at MSU and was a key player in building this Department. Earlier he served as the Associate Chair for Graduate Studies in the Department of Electrical and Computer Engineering from 2012-2015, and the Associate Chair for Research in ECE from 2019-2022. He has authored/co-authored around 450 journal and conference papers and presented a number of invited talks. His research interests include all aspects of computational electromagnetics (frequency and time domain integral equation based methods, multi-scale fast multipole methods, fast transient methods, higher order finite element and integral equation methods), propagation in complex media, mesoscale electromagnetics, and particle and molecular dynamics as applied to multiphysics and multiscale problems. He was an Associate Editor for IEEE Antennas and Wireless Propagation Letters (AWPL), IEEE Transactions on Antennas and Propagation, and Topical Editor for Journal of Optical Society of America: A. He is a full member of the USNC-URSI Commission B. He is Fellow of IEEE (class 2010), elected for his contributions to time and frequency domain computational electromagnetics. He has also been awarded the Withrow Distinguished Junior scholar (in 2003), Withrow Distinguished Senior scholar (in 2010), the Withrow teaching award (in 2007), and the Beal Outstanding Faculty award (2014).
\end{IEEEbiography}

\end{document}